Ion Patrascu | Florentin Smarandache

# Complements to Classic Topics of Circles Geometry

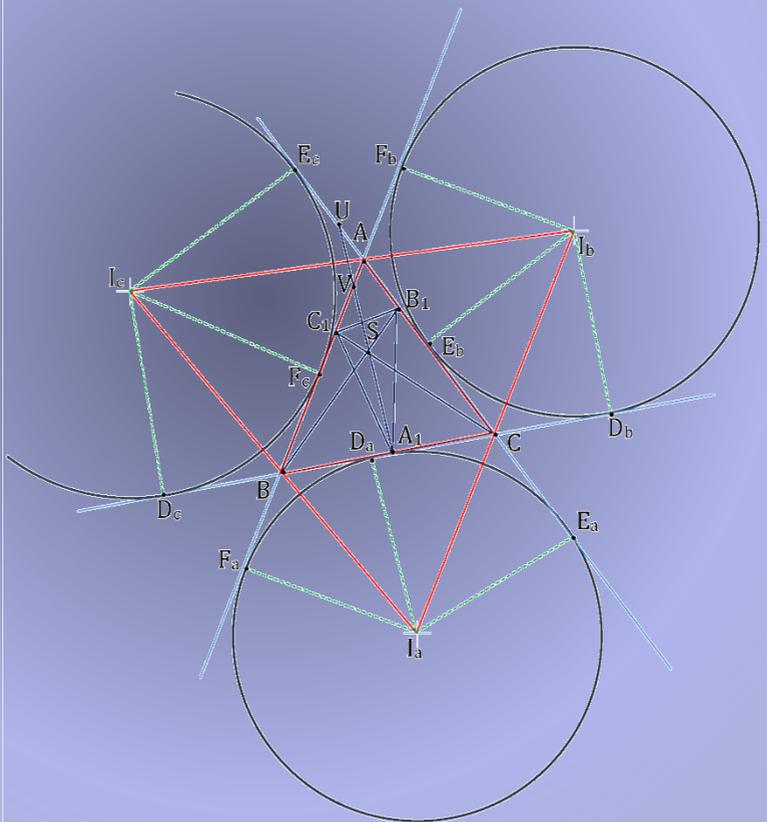

Ion Patrascu | Florentin Smarandache

# Complements to Classic Topics of Circles Geometry

*In the memory of the first author's father Mihail Patrascu and the second author's mother Maria (Marioara) Smarandache, recently passed to eternity...*

Ion Patrascu | Florentin Smarandache

# Complements to Classic Topics of Circles Geometry

# Contents

# **Introductory Note**

We approach several themes of classical geometry of the circle and complete them with some original results, showing that not everything in traditional math is revealed, and that it still has an open character.

The topics were chosen according to authors' aspiration and attraction, as a poet writes lyrics about spring according to his emotions.

# Lemoine's Circles

In this article, we get to **Lemoine's circles** in a different manner than the known one.

## 1$^{st}$ Theorem.

Let $ABC$ a triangle and $K$ its simedian center. We take through K the parallel $A_1A_2$ to $BC$, $A_1 \in (AB)$, $A_2 \in (AC)$; through $A_2$ we take the antiparallels $A_2B_1$ to $AB$ in relation to $CA$ and $CB$, $B_1 \in (BC)$; through $B_1$ we take the parallel $B_1B_2$ to $AC$, $B_2 \in AB$; through $B_2$ we take the antiparallels $B_1C_1$ to $BC$, $C_1 \in (AC)$, and through $C_1$ we take the parallel $C_1C_2$ to $AB$, $C_1 \in (BC)$. Then:

    i.    $C_2A_1$ is an antiparallel of $AC$;

    ii.    $B_1B_2 \cap C_1C_2 = \{K\}$;

    iii.    The points $A_1, A_2$ , $B_1$ , $B_2, C_1, C_2$ are concyclical (the first Lemoine's circle).

*Proof.*

    i.    The quadrilateral $BC_2KA$ is a parallelogram, and its center, i.e. the middle of the segment $(C_2A_1)$, belongs to the simedian $BK$; it follows that $C_2A_2$ is an antiparallel to $AC$ (see *Figure 1*).

ii. Let $\{K'\} = A_1A_2 \cap B_1B_2$ , because the quadrilateral $K'B_1CA_2$ is a parallelogram; it follows that $CK'$ is a simedian; also, $CK$ is a simedian, and since $K, K' \in A_1A_2$, it follows that we have $K' = K$.

iii. $B_2C_1$ being an antiparallel to $BC$ and $A_1A_2 \parallel BC$ , it means that $B_2C_1$ is an antiparallel to $A_1A_2$ , so the points $B_2, C_1, A_2, A_1$ are concyclical. From $B_1B_2 \parallel AC$ , $\sphericalangle B_2C_1A \equiv \sphericalangle ABC$ , $\sphericalangle B_1A_2C \equiv \sphericalangle ABC$ we get that the quadrilateral $B_2C_1A_2B_1$ is an isosceles trapezoid, so the points $B_2, C_1, A_2, B_1$ are concyclical. Analogously, it can be shown that the quadrilateral $C_2B_1A_2A_1$ is an isosceles trapezoid, therefore the points $C_2, B_1, A_2, A_1$ are concyclical.

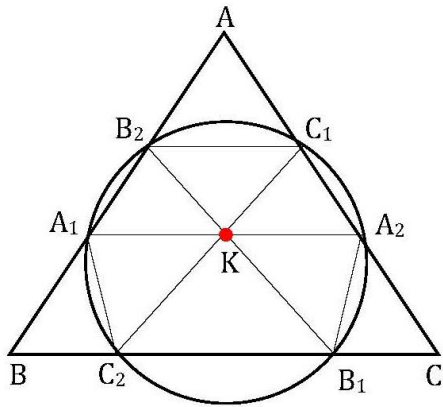

*Figure 1*

From the previous three quartets of concyclical points, it results the concyclicity of the points belonging to the first Lemoine's circle.

## 2$^{nd}$ Theorem.

In the scalene triangle $ABC$, let K be the simedian center. We take from $K$ the antiparallel $A_1A_2$ to $BC$; $A_1 \in AB, A_2 \in AC$; through $A_2$ we build $A_2B_1 \parallel AB$; $B_1 \in (BC)$, then through $B_1$ we build $B_1B_2$ the antiparallel to $AC$, $B_2 \in (AB)$, and through $B_2$ we build $B_2C_1 \parallel BC$, $C_1 \in AC$, and, finally, through $C_1$ we take the antiparallel $C_1C_2$ to $AB$, $C_2 \in (BC)$.

Then:

i. $C_2A_1 \parallel AC$;

ii. $B_1B_2 \cap C_1C_2 = \{K\}$;

iii. The points $A_1, A_2, B_1, B_2, C_1, C_2$ are concyclical (the second Lemoine's circle).

*Proof.*

i. Let $\{K'\} = A_1A_2 \cap B_1B_2$, having $\sphericalangle AA_1A_2 = \sphericalangle ACB$ and $\sphericalangle BB_1B_2 \equiv \sphericalangle BAC$ because $A_1A_2$ și $B_1B_2$ are antiparallels to $BC$, $AC$, respectively, it follows that $\sphericalangle K'A_1B_2 \equiv \sphericalangle K'B_2A_1$, so $K'A_1 = K'B_2$; having $A_1B_2 \parallel B_1A_2$ as well, it follows that also $K'A_2 = K'B_1$, so $A_1A_2 = B_1B_2$. Because $C_1C_2$ and $B_1B_2$ are antiparallels to $AB$ and $AC$, we have $K"C_2 = K"B_1$; we noted $\{K"\} = B_1B_2 \cap C_1C_2$; since $C_1B_2 \parallel B_1C_2$, we have that the triangle $K"C_1B_2$ is also isosceles, therefore $K"C_1 = C_1B_2$, and we get that $B_1B_2 = C_1C_2$. Let $\{K'''\} = A_1A_2 \cap C_1C_2$ ; since $A_1A_2$ and $C_1C_2$ are antiparallels to $BC$ and $AB$, we get that the triangle $K'''A_2C_1$ is isosceles, so $K'''A_2 = K'''C_1$, but $A_1A_2 = C_1C_2$ implies that $K'''C_2 = K'''A_1$ , then $\sphericalangle K'''A_1C_2 \equiv \sphericalangle K'''A_2C_1$ and, accordingly, $C_2A_1 \parallel AC$.

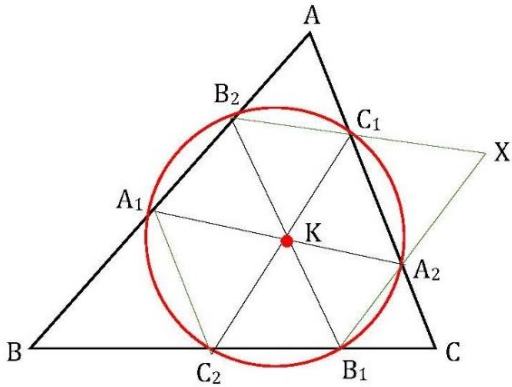

*Figure 2*

ii.      We noted $\{K'\} = A_1A_2 \cap B_1B_2$ ; let $\{X\} = B_2C_1 \cap B_1A_2$ ; obviously, $BB_1XB_2$ is a parallelogram; if $K_0$ is the middle of $(B_1B_2)$, then $BK_0$ is a simedian, since $B_1B_2$ is an antiparallel to $AC$, and the middle of the antiparallels of $AC$ are situated on the simedian $BK$. If $K_0 \neq K$, then $K_0 K \parallel A_1 B_2$ (because $A_1 A_2 = B_1 B_2$ and $B_1 A_2 \parallel A_1 B_2$ ), on the other hand, $B, K_0, K$ are collinear (they belong to the simedian $BK$ ), therefore $K_0 K$ intersects $AB$ in $B$, which is absurd, so $K_0 = K$, and, accordingly, $B_1 B_2 \cap A_1 A_2 = \{K\}$. Analogously, we prove that $C_1 C_2 \cap A_1 A_2 = \{K\}$, so $B_1 B_2 \cap C_1 C_2 = \{K\}$.

iii. K is the middle of the congruent antiparalells $A_1 A_2$ , $B_1 B_2$ , $C_1 C_2$ , so $KA_1 = KA_2 = KB_1 = KB_2 = KC_1 = KC_2$ . The simedian center $K$ is the center of the second Lemoine's circle.

### *Remark.*

The center of the first Lemoine's circle is the middle of the segment $[OK]$, where $O$ is the center of the circle circumscribed to the triangle $ABC$. Indeed, the perpendiculars taken from $A, B, C$ on the antiparallels $B_2 C_1$ , $A_1 C_2$ , $B_1 A_2$ respectively pass through O, the center of the circumscribed circle (the antiparallels have the directions of the tangents taken to the circumscribed circle in $A, B, C$). The mediatrix of the segment $B_2 C_1$ pass though the middle of $B_2 C_1$ , which coincides with the middle of $AK$, so is the middle line in the triangle $AKO$ passing through the middle of $(OK)$. Analogously, it follows that the mediatrix of $A_1 C_2$ pass through the middle $L_1$ of $[OK]$.

# References.

# Lemoine's Circles Radius Calculus

For the calculus of the **first Lemoine's circle**, we will first prove:

## 1$^{st}$ Theorem

(E. Lemoine – 1873)

The first Lemoine's circle divides the sides of a triangle in segments proportional to the squares of the triangle's sides.

Each extreme segment is proportional to the corresponding adjacent side, and the chord-segment in the Lemoine's circle is proportional to the square of the side that contains it.

*Proof.*

We will prove that $\frac{BC_2}{c^2} = \frac{C_2B_1}{a^2} = \frac{B_1C}{b^2}$ .

In figure 1, $K$ represents the symmedian center of the triangle $ABC$, and $A_1A_2$ ; $B_1B_2$ ; $C_1C_2$ represent Lemoine parallels.

The triangles $BC_2A_1$ ; $CB_1A_2$ and $KC_2A_1$ have heights relative to the sides $BC_2$ ; $B_1C$ and $C_2B_1$ equal ($A_1A_2 \parallel BC$).

Hence:

$$\frac{Area_\Delta BA_1C_2}{BC_2} = \frac{Area_\Delta KC_2A_1}{C_2B_1} = \frac{Area_\Delta CB_1A_2}{B_1C} \ . \tag{1}$$

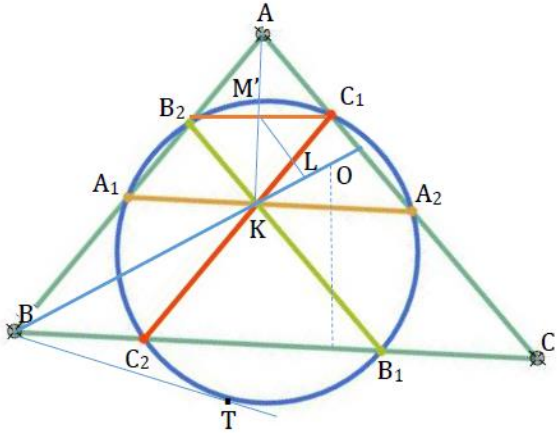

*Figure 1*

On the other hand: $A_1C_2$ and $B_1A_2$ being antiparallels with respect to $AC$ and $AB$, it follows that $\Delta BC_2A_1 \sim \Delta BAC$ and $\Delta CB_1A_2 \sim \Delta CAB$, likewise $KC_2 \parallel AC$ implies: $\Delta KC_2B_1 \sim \Delta ABC$.

We obtain:

$$\frac{Area_\Delta BC_2A_1}{Area_\Delta ABC} = \frac{BC_2^2}{c^2} \ ;$$

$$\frac{Area_\Delta KC_2B_1}{Area_\Delta ABC} = \frac{C_2B_1^2}{a^2} \ ;$$

$$\frac{Area_\Delta CB_1A_2}{Area_\Delta ABC} = \frac{CB_1^2}{b^2} \ . \tag{2}$$

If we denote $Area_\Delta ABC = S$, we obtain from the relations (1) and (2) that:

$$\frac{BC_2}{c^2} = \frac{C_2B_1}{a^2} = \frac{B_1C}{b^2} \ .$$

*Consequences.*

1. According to the 1$^{st}$ Theorem, we find that:
$BC_2 = \frac{ac^2}{a^2+b^2+c^2}$; $B_1C = \frac{ab^2}{a^2+b^2+c^2}$; $B_1C_2 = \frac{a^3}{a^2+b^2+c^2}$ .

2. We also find that:
$\frac{B_1C_2}{a^3} = \frac{A_2C_1}{b^3} = \frac{A_1B_2}{c^3}$ ,

meaning that:

*"The chords determined by the first Lemoine's circle on the triangle's sides are proportional to the cubes of the sides."*

Due to this property, the first Lemoine's circle is known in England by the name of *triplicate ratio circle*.

# 1$^{st}$ **Proposition**.

The radius of the first Lemoine's circle, $R_{L_1}$ is given by the formula:

$$R_{L_1}^2 = \frac{1}{4} \cdot \frac{R^2(a^2+b^2+c^2) + a^2b^2c^2}{(a^2+b^2+c^2)^2}, \tag{3}$$

where $R$ represents the radius of the circle inscribed in the triangle $ABC$.

*Proof.*

Let $L$ be the center of the first Lemoine's circle that is known to represent the middle of the segment $(OK)$ – $O$ being the center of the circle inscribed in the triangle $ABC$.

Considering C1, we obtain $BB_1 = \frac{a\,(c^2+a^2)}{a^2+b^2+c^2}$ .

Taking into account the power of point $B$ in relation to the first Lemoine's circle, we have:

$$BC_2 \cdot BB_1 = BT^2 - LT^2,$$

($BT$ is the tangent traced from $B$ to the first Lemoine's circle, see *Figure 1*).

Hence: $R_{L_1}^2 = BL^2 - BC_2 \cdot BB_1.$     (4)

The median theorem in triangle $BOK$ implies that:

$$BL^2 = \frac{2 \cdot \left(BK^2 + BO^2\right) - OK^2}{4} \ .$$

It is known that $K = \frac{\left(a^2+c^2\right) \cdot S_b}{a^2+b^2+c^2}$ ; $S_b = \frac{2ac \cdot m_b}{a^2+c^2}$ , where $S_b$ and $m_b$ are the lengths of the symmedian and the median from $B$, and $OK^2 = R^2 - \frac{3a^2b^2c^2}{(a^2+b^2+c^2)}$ , see (3).

Consequently: $BK^2 = \frac{2a^2c^2\left(a^2+c^2\right) - a^2b^2c^2}{(a^2+b^2+c^2)^2}$ , and

$$4BL^2 = R^2 + \frac{4a^2c^2\left(a^2+c^2\right) + a^2b^2c^2}{(a^2+b^2+c^2)^2} \ .$$

As: $BC_2 \cdot BB_1 = \frac{a^2c^2\left(a^2+c^2\right)}{(a^2+b^2+c^2)^2}$ , by replacing in (4), we obtain formula (3).

## 2$^{\text{nd}}$ Proposition.

The radius of the second Lemoine's circle, $R_{L_2}$, is given by the formula:

$$R_{L_2} = \frac{abc}{a^2+b^2+c^2} \ . \tag{5}$$

*Proof.*

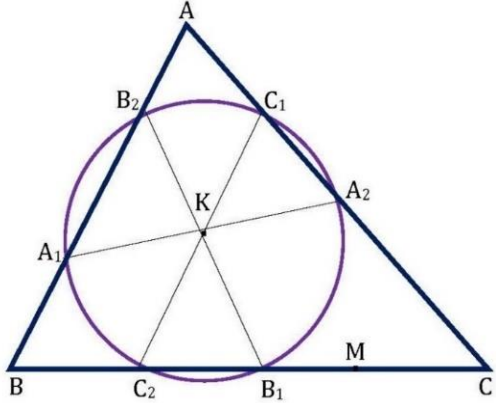

*Figure 2*

In *Figure 2* above, $A_1A_2$; $B_1B_2$; $C_1C_2$ are Lemoine antiparallels traced through symmedian center $K$ that is the center of the second Lemoine's circle, thence:

$R_{L_2} = KA_1 = KA_2$.

If we note with $S$ and $M$ the feet of the symmedian and the median from $A$, it is known that:

$\frac{AK}{KS} = \frac{b^2+c^2}{a^2}$ .

From the similarity of triangles $AA_2A_1$ and $ABC$, we have: $\frac{A_1A_2}{BC} = \frac{AK}{AM}$ .

But: $\frac{AK}{AS} = \frac{b^2+c^2}{a^2+b^2+c^2}$ and $AS = \frac{2bc}{b^2+c^2} \cdot m_a$.

$A_1A_2 = 2R_{L_2}$, $BC = a$, therefore:

$R_{L_2} = \frac{AK \cdot a}{2m_a}$ ,

and as $AK = \frac{2bc \cdot m_a}{a^2+b^2+c^2}$ , formula (5) is a consequence.

*Remarks.*

1. If we use $tg\omega = \frac{4S}{a^2+b^2+c^2}$, $\omega$ being the Brocard's angle (see [2]), we obtain: $R_{L_2} = R \cdot tg\omega$.

2. If, in *Figure 1*, we denote with $M_1$ the middle of the antiparallel $B_2C_1$, which is equal to $R_{L_2}$ (due to their similarity), we thus find from the rectangular triangle $LM_1C_1$ that:

$LC_1^2 = LM_1^2 + M_1C_1^2$ , but $LM_1^2 = \frac{1}{4}a^2$ and $M_1C_2 = \frac{1}{2}R_{L_2}$; it follows that:

$R_{L_1}^2 = \frac{1}{4}\left(R^2 + R_{L_2}^2\right) = \frac{R^2}{4}(1 + tg^2\omega)$.

We obtain:

$R_{L_1} = \frac{R}{2} \cdot \sqrt{1 + tg^2\omega}$ .

# 3[rd] **Proposition.**

The chords determined by the sides of the triangle in the second Lemoine's circle are respectively proportional to the opposing angles cosines.

*Proof.*

$KC_2B_1$ is an isosceles triangle, $\sphericalangle KC_2B_1 = \sphericalangle KB_1C_2 = \sphericalangle A$; as $KC_2 = R_{L_2}$ we have that $\cos A = \frac{C_2B_1}{2R_{L_2}}$, deci $\frac{C_2B_1}{\cos A} = 2R_{L_2}$, similarly: $\frac{A_2C_1}{\cos B} = \frac{B_2A_1}{\cos C} = 2R_{L_2}$.

*Remark.*

Due to this property of the Lemoine's second circle, in England this circle is known as the *cosine circle.*

# References.

# Radical Axis of Lemoine's Circles

In this article, we emphasize the **radical axis of the Lemoine's circles**.

For the start, let us remind:

## 1st Theorem.

The parallels taken through the simmedian center $K$ of a triangle to the sides of the triangle determine on them six concyclic points (the first Lemoine's circle).

## 2nd Theorem.

The antiparallels taken through the triangle's simmedian center to the sides of a triangle determine six concyclic points (the second Lemoine's circle).

### *1st Remark.*

If $ABC$ is a scalene triangle and $K$ is its simmedian center, then $L$, the center of the first Lemoine's circle, is the middle of the segment $[OK]$, where $O$ is the center of the circumscribed circle, and the center of the second Lemoine's circle is $K$. It follows that the radical axis of Lemoine's circles is perpendicular on the line of the centers $LK$, therefore on the line $OK$.

## 1$^{st}$ Proposition.

The radical axis of Lemoine's circles is perpendicular on the line $OK$ raised in the simmedian center $K$.

### Proof.

Let $A_1A_2$ be the antiparallel to $BC$ taken through $K$, then $KA_1$ is the radius $R_{L_2}$ of the second Lemoine's circle; we have:

$R_{L_2} = \frac{abc}{a^2+b^2+c^2}$ .

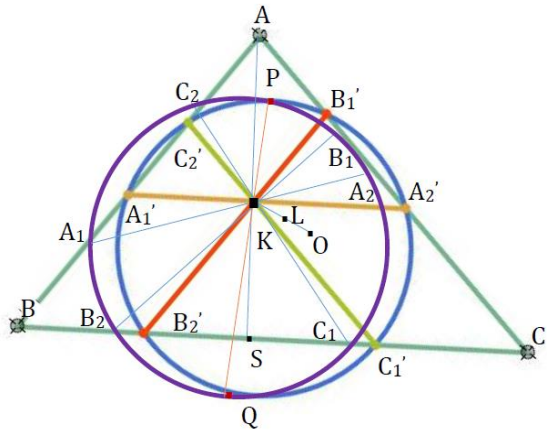

*Figure 1*

Let $A_1' A_2'$ be the Lemoine's parallel taken to $BC$; we evaluate the power of $K$ towards the first Lemoine's circle. We have:

$$\overrightarrow{KA_1'} \cdot \overrightarrow{KA_2'} = LK^2 - R_{L_1}^2. \tag{1}$$

Let $S$ be the simmedian leg from $A$; it follows that:

$\frac{KA_1'}{BS} = \frac{AK}{AS} - \frac{KA_2'}{SC}$ .

We obtain:

$KA_1' = BS \cdot \frac{AK}{AS}$ and $KA_2' = SC \cdot \frac{AK}{AS}$ ,

but $\frac{BS}{SC} = \frac{c^2}{b^2}$ and $\frac{AK}{AS} = \frac{b^2 + c^2}{a^2 + b^2 + c^2}$ .

Therefore:

$$\overrightarrow{KA_1'} \cdot \overrightarrow{KA_2'} = -BS \cdot SC \cdot \left(\frac{AK}{AS}\right)^2 = \frac{-a^2 b^2 c^2}{(b^2 + c^2)^2};$$

$$\frac{(b^2 + c^2)^2}{(a^2 + b^2 + c^2)^2} = -R_{L_2}^2 . \tag{2}$$

We draw the perpendicular in $K$ on the line $LK$ and denote by $P$ and $Q$ its intersection to the first Lemoine's circle.

We have $\overrightarrow{KP} \cdot \overrightarrow{KQ} = -R_{L_2}^2$ ; by the other hand, $KP = KQ$ ($PQ$ is a chord which is perpendicular to the diameter passing through $K$).

It follows that $KP = KQ = R_{L_2}$, so $P$ and $Q$ are situated on the second Lemoine's circle.

Because $PQ$ is a chord which is common to the Lemoine's circles, it obviously follows that $PQ$ is the radical axis.

*Comment.*

Equalizing (1) and (2), or by the Pythagorean theorem in the triangle $PKL$, we can calculate $R_{L_1}$.

It is known that: $OK^2 = R^2 - \frac{3a^2b^2c^2}{(a^2+b^2+c^2)^2}$ , and since $LK = \frac{1}{2}OK$, we find that:

$$R_{L_1}^2 = \frac{1}{4} \cdot \left[ R^2 + \frac{a^2b^2c^2}{(a^2+b^2+c^2)^2} \right].$$

*$2^{nd}$ Remark.*

The 1[st] Proposition, ref. the radical axis of the Lemoine's circles, is a particular case of the following *Proposition*, which we leave to the reader to prove.

## 2[nd] Proposition.

If $\mathcal{C}(O_1, R_1)$ şi $\mathcal{C}(O_2, R_2)$ are two circles such as the power of center $O_1$ towards $\mathcal{C}(O_2, R_2)$ is $-R_1^2$, then the radical axis of the circles is the perpendicular in $O_1$ on the line of centers $O_1O_2$.

## References.

# Generating Lemoine's circles

In this paper, we generalize the theorem relative to the first Lemoine's circle and thereby highlight **a method to build Lemoine's circles**.

Firstly, we review some notions and results.

## 1st Definition.

It is called a simedian of a triangle the symmetric of a median of the triangle with respect to the internal bisector of the triangle that has in common with the median the peak of the triangle.

## 1st Proposition.

In the triangle $ABC$, the cevian $AS$, $S \in (BC)$, is a simedian if and only if $\frac{SB}{SC} = \left(\frac{AB}{AC}\right)^2$. For *Proof*, see [2].

## 2nd Definition.

It is called a simedian center of a triangle (or Lemoine's point) the intersection of triangle's simedians.

# 1ˢᵗ Theorem.

The parallels to the sides of a triangle taken through the simedian center intersect the triangle's sides in six concyclic points (the first Lemoine's circle - 1873).

A *Proof* of this theorem can be found in [2].

# 3ʳᵈ Definition.

We assert that in a scalene triangle $ABC$ the line $MN$, where $M \in AB$ and $N \in AC$, is an antiparallel to $BC$ if $\sphericalangle MNA \equiv \sphericalangle ABC$.

# 1ˢᵗ Lemma.

In the triangle $ABC$, let $AS$ be a simedian, $S \in (BC)$. If $P$ is the middle of the segment $(MN)$, having $M \in (AB)$ and $N \in (AC)$, belonging to the simedian $AS$, then $MN$ and $BC$ are antiparallels.

## *Proof.*

We draw through $M$ and $N$, $MT \parallel AC$ and $NR \parallel AB$, $R,T \in (BC)$, see *Figure 1*. Let $\{Q\} = MT \cap NR$; since $MP = PN$ and $AMQN$ is a parallelogram, it follows that $Q \in AS$.

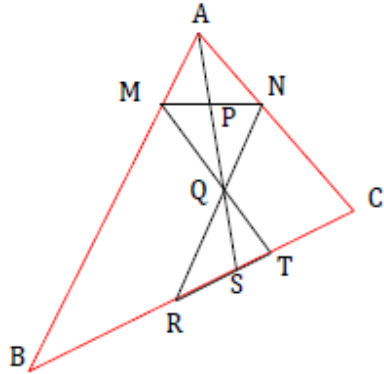

*Figure 1.*

Thales's Theorem provides the relations:

$$\frac{AN}{AC} = \frac{BR}{BC};\tag{1}$$

$$\frac{AB}{AM} = \frac{BC}{CT}.\tag{2}$$

From (1) and (2), by multiplication, we obtain:

$$\frac{AN}{AM} \cdot \frac{AB}{AC} = \frac{BR}{TC}.\tag{3}$$

Using again Thales's Theorem, we obtain:

$$\frac{BR}{BS} = \frac{AQ}{AS},\tag{4}$$

$$\frac{TC}{SC} = \frac{AQ}{AS}.\tag{5}$$

From these relations, we get

$$\frac{BR}{BS} = \frac{TC}{SC},\tag{6}$$

or

$$\frac{BS}{SC} = \frac{BR}{TC}.\tag{7}$$

In view of *Proposition 1*, the relations (7) and (3) lead to $\frac{AN}{AB} = \frac{AB}{AC}$, which shows that $\Delta AMN \sim \Delta ACB$, so $\measuredangle AMN \equiv \measuredangle ABC$.

Therefore, $MN$ and $BC$ are antiparallels in relation to $AB$ and $AC$.

### $Remark.$

1. The reciprocal of *Lemma 1* is also valid, meaning that if $P$ is the middle of the antiparallel $MN$ to $BC$, then $P$ belongs to the simedian from $A$.

## $2^{nd}$ Theorem.

(Generalization of the *$1^{st}$ Theorem*)

Let $ABC$ be a scalene triangle and $K$ its simedian center. We take $M \in AK$ and draw $MN \parallel AB, MP \parallel AC$, where $N \in BK$, $P \in CK$. Then:

    i.      $NP \parallel BC$;

    ii.     $MN, NP$ and $MP$ intersect the sides of triangle $ABC$ in six concyclic points.

### $Proof.$

In triangle $ABC$, let $AA_1$, $BB_1, CC_1$ the simedians concurrent in $K$ (see *Figure 2*).

We have from Thales' Theorem that:

$$\frac{AM}{MK} = \frac{BN}{NK}; \tag{1}$$

$$\frac{AM}{MK} = \frac{CP}{PK} . \tag{2}$$

From relations (1) and (2), it follows that

$$\frac{BN}{NK} = \frac{CP}{PK}, \tag{3}$$

which shows that $NP \parallel BC$.

Let $R, S, V, W, U, T$ be the intersection points of the parallels $MN, MP, NP$ of the sides of the triangles to the other sides.

Obviously, by construction, the quadrilaterals $ASMW$; $CUPV$; $BRNT$ are parallelograms.

The middle of the diagonal $WS$ falls on $AM$, so on the simedian $AK$, and from $1^{st}$ *Lemma* we get that $WS$ is an antiparallel to $BC$.

Since $TU \parallel BC$, it follows that $WS$ and $TU$ are antiparallels, therefore the points $W, S, U, T$ are concyclic (4).

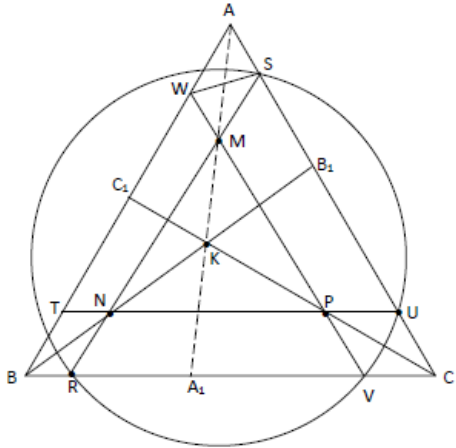

*Figure 2.*

Analogously, we show that the points $U, V, R, S$ are concyclic (5). From $WS$ and $BC$ antiparallels, $UV$ and $AB$ antiparallels, we have that $\sphericalangle WSA \equiv \sphericalangle ABC$ and $\sphericalangle VUC \equiv \sphericalangle ABC$, therefore: $\sphericalangle WSA \equiv \sphericalangle VUC$, and since $VW \parallel AC$, it follows that the trapeze $WSUV$ is isosceles, therefore the points $W, S, U, V$ are concyclic (6).

The relations (4), (5), (6) drive to the concyclicality of the points $R, U, V, S, W, T$, and the theorem is proved.

### *Further Remarks.*

2. For any point $M$ found on the simedian $AA_1$, by performing the constructions from hypothesis, we get a circumscribed circle of the 6 points of intersection of the parallels taken to the sides of triangle.

3. The $2^{nd}$ *Theorem* generalizes the $1^{st}$ *Theorem* because we get the second in the case the parallels are taken to the sides through the simedian center $k$.

4. We get a circle built as in $2^{nd}$ *Theorem* from the first Lemoine's circle by homothety of pole $k$ and of ratio $\lambda \in \mathbb{R}$.

5. The centers of Lemoine's circles built as above belong to the line $OK$, where $O$ is the center of the circle circumscribed to the triangle $ABC$.

## References.

# The Radical Circle of Ex-Inscribed Circles of a Triangle

In this article, we prove several theorems about the **radical center** and the **radical circle of ex-inscribed circles of a triangle** and calculate the **radius of the circle from vectorial considerations**.

## 1$^{\text{st}}$ Theorem.

The radical center of the ex-inscribed circles of the triangle $ABC$ is the Spiecker's point of the triangle (the center of the circle inscribed in the median triangle of the triangle $ABC$).

### *Proof.*

We refer in the following to the notation in *Figure 1*. Let $I_a, I_b, I_c$ be the centers of the ex-inscribed circles of a triangle (the intersections of two external bisectors with the internal bisector of the other angle). Using tangents property taken from a point to a circle to be congruent, we calculate and find that:

$$AF_a = AE_a = BD_b = BF_b = CD_c = CE_c = p,$$
$$BD_c = BF_c = CD_b = CE_b = p - a,$$

$$CE_a = CD_a = AF_c = AE_c = p - b,$$
$$AF_b = AE_b = BF_c = BD_c = p - c.$$

If $A_1$ is the middle of segment $D_cD_b$, it follows that $A_1$ has equal powers to the ex-inscribed circles $(I_b)$ and $(I_c)$. Of the previously set, we obtain that $A_1$ is the middle of the side $BC$.

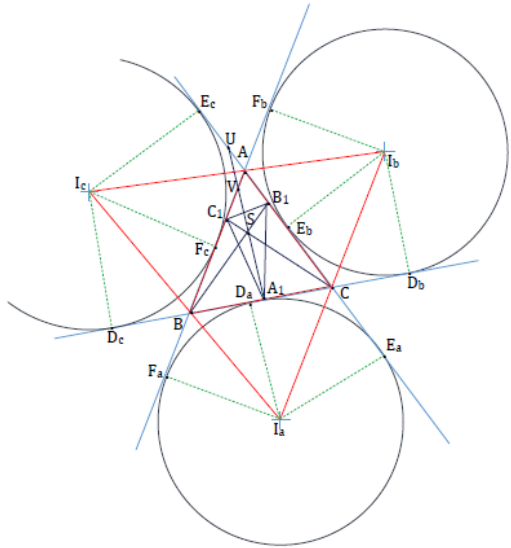

*Figure 1.*

Also, the middles of the segments $E_bE_c$ and $F_bF_c$, which we denote $U$ and $V$, have equal powers to the circles $(I_b)$ and $(I_c)$.

The radical axis of the circles $(I_b)$, $(I_c)$ will include the points $A_1, U, V$.

Because $AE_b = AF_b$ and $AE_c = AF_c$, it follows that $AU = AY$ and we find that $\sphericalangle AUV = \frac{1}{2}\sphericalangle A$, therefore the radical axis of the ex-inscribed circles $(F_b)$ and $(I_c)$ is the parallel taken through the middle $A_1$ of the side $BC$ to the bisector of the angle $BAC$.

Denoting $B_1$ and $C_1$ the middles of the sides $AC$, $AB$, respectively, we find that the radical center of the ex-inscribed circles is the center of the circle inscribed in the median triangle $A_1B_1C_1$ of the triangle $ABC$.

This point, denoted $S$, is the Spiecker's point of the triangle ABC.

## 2$^{nd}$ Theorem.

The radical center of the inscribed circle $(I)$ and of the $B$ −ex-inscribed and $C$ −ex-inscribed circles of the triangle $ABC$ is the center of the $A_1$ − ex-inscribed circle of the median triangle $A_1B_1C_1$, corresponding to the triangle $ABC$).

### *Proof.*

If $E$ is the contact of the inscribed circle with $AC$ and $E_b$ is the contact of the $B$ −ex-inscribed circle with $AC$ , it is known that these points are isotomic, therefore the middle of the segment $EE_b$ is the middle of the side $AC$, which is $B_1$.

This point has equal powers to the inscribed circle $(I)$ and to the $B$ −ex-inscribed circle $(I_b)$, so it belongs to their radical axis.

Analogously, $C_1$ is on the radical axis of the circles $(I)$ and $(I_c)$.

The radical axis of the circles $(I)$, $(I_b)$ is the perpendicular taken from $B_1$ to the bisector $II_b$.

This bisector is parallel with the internal bisector of the angle $A_1B_1C_1$ , therefore the perpendicular in $B_1$ on $II_b$ is the external bisector of the angle $A_1B_1C_1$ from the median triangle.

Analogously, it follows that the radical axis of the circles $(I)$, $(I_c)$ is the external bisector of the angle $A_1C_1B_1$ from the median triangle.

Because the bisectors intersect in the center of the circle $A_1$ -ex-inscribed to the median triangle $A_1B_1C_1$, this point $S_a$ is the center of the radical center of the circles $(I)$, $(I_b)$, $(I_c)$.

### Remark.

The theorem for the circles $(I)$, $(I_a)$, $(I_b)$ and $(I)$, $(I_a)$, $(I_c)$ can be proved analogously, obtaining the points $S_c$ and $S_b$.

## 3$^{\text{rd}}$ Theorem.

The radical circle's radius of the circles ex-inscribed to the triangle $ABC$ is given by the formula: $\frac{1}{2}\sqrt{r^2 + p^2}$, where $r$ is the radius of the inscribed circle.

*Proof.*

The position vector of the circle $I$ of the inscribed circle in the triangle ABC is:

$$\overrightarrow{PI} = \frac{1}{2p}(a\overrightarrow{PA} + b\overrightarrow{PB} + c\overrightarrow{PC}).$$

Spiecker's point $S$ is the center of radical circle of ex-inscribed circle and is the center of the inscribed circle in the median triangle $A_1B_1C_1$, therefore:

$$\overrightarrow{PS} = \frac{1}{p}\left(\frac{1}{2}a\overrightarrow{PA_1} + \frac{1}{2}b\overrightarrow{PB_1} + \frac{1}{2}c\overrightarrow{PC_1}\right).$$

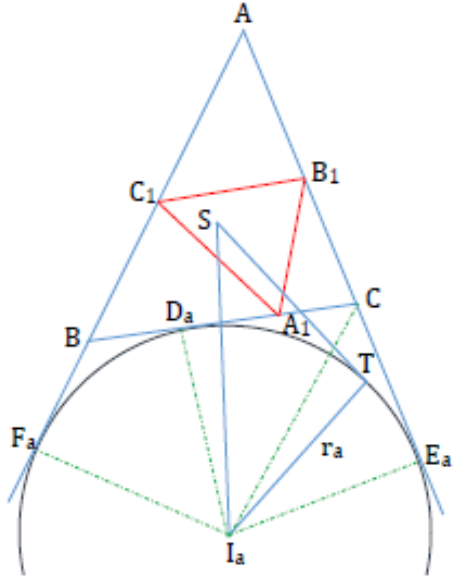

*Figure 2.*

We denote by $T$ the contact point with the $A$-ex-inscribed circle of the tangent taken from $S$ to this circle (see *Figure 2*).

The radical circle's radius is given by:

$$ST = \sqrt{SI_a^2 - I_a^2}$$

$$\overrightarrow{I_aS} = \frac{1}{2p}\left(a\overrightarrow{I_aA_1} + b\overrightarrow{I_aB_1} + c\overrightarrow{I_aC_1}\right).$$

We evaluate the product of the scales $\overrightarrow{I_aS} \cdot \overrightarrow{I_aS}$ ; we have:

$$I_aS^2 = \frac{1}{4p^2}\left(a^2I_aA_1^2 + b^2I_aB_1^2 + c^2I_aC_1^2 + 2ab\overrightarrow{I_aA_1} \cdot \overrightarrow{I_aB_1} + 2bc\overrightarrow{I_aB_1} \cdot \overrightarrow{I_aC_1} + 2ac\overrightarrow{I_aA_1} \cdot \overrightarrow{I_aC_1}\right).$$

From the law of cosines applied in the triangle $I_aA_1B_1$, we find that:

$$2\overrightarrow{I_aA_1} \cdot \overrightarrow{I_aB_1} = I_aA_1^2 + I_aB_1^2 - \frac{1}{4}c^2, \text{ therefore:}$$

$$2ab\overrightarrow{I_aA_1} \cdot \overrightarrow{I_aB_1} = ab(I_aA_1^2 + I_aB_1^2 - \frac{1}{4}abc^2.$$

Analogously, we obtain:

$$2bc\overrightarrow{I_aB_1} \cdot \overrightarrow{I_aC_1} = bc(I_aB_1^2 + I_aC_1^2 - \frac{1}{4}a^2bc,$$

$$2ac\overrightarrow{I_aA_1} \cdot \overrightarrow{I_aC_1} = ac(I_aA_1^2 + I_aC_1^2 - \frac{1}{4}ab^2c.$$

$$I_aS^2 = \frac{1}{4p^2}\Big[(a^2 + ab + ac)I_aA_1^2 + (b^2 + ab + bc)I_aB_1^2 + (c^2 + bc + ac)I_aC_1^2 - \frac{abc}{4}(a + b + c)\Big],$$

$$I_aS^2 = \frac{1}{4p^2}\Big[2p(aI_aA_1^2 + bI_aB_1^2 + cI_aC_1^2) - 2RS_p\Big],$$

$$I_aS^2 = \frac{1}{2p}(aI_aA_1^2 + bI_aB_1^2 + cI_aC_1^2) - \frac{1}{2}Rr.$$

From the right triangle $I_aD_aA_1$, we have that:

$$I_aA_1^2 = r_a^2 + A_1D_a^2 = r_a^2 + \left[\frac{a}{2} - (p - c)\right]^2 =$$
$$= r_a^2 + \frac{(c-b)^2}{4}\,.$$

From the right triangles $I_aE_aB_1$ și $I_aF_aC_1$, we find:

$$I_a B_1^2 = r_a^2 + B_1 E_a^2 = r_a^2 + \left[\frac{b}{2} - (p - b)\right]^2 =$$
$$= r_a^2 + \frac{1}{4}(a + c)^2,$$

$$I_a C_1^2 = r_a^2 + \frac{1}{4}(a + b)^2.$$

Evaluating $aI_a A_1^2 + bI_a B_1^2 + cI_a C_1^2$, we obtain:

$$aI_a A_1^2 + bI_a B_1^2 + cI_a C_1^2 =$$
$$= 2pr_a^2 + \frac{1}{2}p(ab + ac + bc) - \frac{1}{4}abc.$$

But:

$$ab + ac + bc = r^2 + p^2 + 4Rr.$$

It follows that:

$$\frac{1}{2p}[aI_a A_1^2 + bI_a B_1^2 + cI_a C_1^2] = r_a^2 + \frac{1}{4}(r^2 + p^2) + \frac{1}{2}Rr$$

and

$$I_a S^2 = r_a^2 + \frac{1}{4}(r^2 + p^2).$$

Then, we obtain:

$$ST = \frac{1}{2}\sqrt{r^2 + p^2}.$$

## References.

# The Polars of a Radical Center

In [1] , the late mathematician **Cezar Cosnita**, using the barycenter coordinates, proves two theorems which are the subject of this article.

In a remark made after proving the first theorem, C. Cosnita suggests an elementary proof by employing the **concept of polar**.

In the following, we prove the theorems based on the indicated path, and state that the second theorem is a particular case of the former. Also, we highlight other particular cases of these theorems.

## 1ˢᵗ Theorem.

Let $ABC$ be a given triangle; through the pairs of points $(B, C)$, $(C, A)$ and $(A, B)$ we take three circles such that their radical center is on the outside.

The polar lines of the radical center of these circles in relation to each of them cut the sides $BC$, $CA$ and $AB$ respectively in three collinear points.

*Proof.*

We denote by $D, P, F$ the second point of intersection of the pairs of circles passing through $(A, B)$ and $(B, C)$; $(B, C)$ and $(A, C)$, $(B, C)$ and $(A, B)$ respectively (see *Figure 1*).

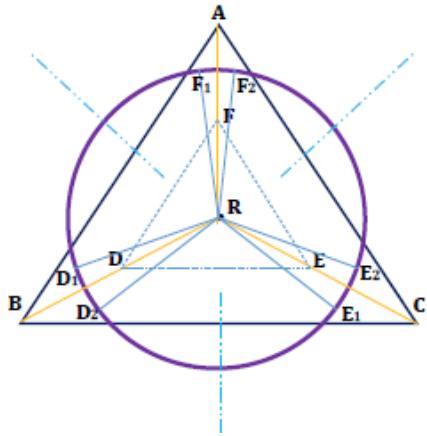

*Figure 1*

Let $R$ be the radical center of those circles. In fact, $\{R\} = AF \cap BD \cap CE$.

We take from $R$ the tangents $RD_1 = RD_2$ to the circle $(B, C)$, $RE_1 = RE_2$ to the circle $(A, C)$ and $RF_1 = RF_2$ to the circle passing through $(A, B)$. Actually, we build the radical circle $\mathcal{C}(R, RD_1)$ of the given circles.

The polar lines of $R$ to these circles are the lines $D_1 D_2$, $E_1 E_2$, $F_1 F_2$. These three lines cut $BC, AC$ and $AB$ in the points $X, Y$ and $Z$, and these lines are respectively the polar lines of $R$ in respect to the circles passing through $(B, C), (C, A)$ and $(A, B)$. The polar lines are the radical axis of the radical circle with each of the circles passing through $(B, C), (C, A), (A, B)$, respectively. The points belong to the radical axis having equal powers to those circles, thereby $XD_1 \cdot XD_2 = XC \cdot XB$.

This relationship shows that the point $X$ has equal powers relative to the radical circle and to the circle circumscribed to the triangle $ABC$; analogically, the point $Y$ has equal powers relative to the radical circle and to the circle circumscribed to the triangle $ABC$; and, likewise, the point $Z$ has equal powers relative to the radical circle and to the circle circumscribed to the triangle $ABC$.

Because the locus of the points having equal powers to two circles is generally a line, i.e. their radical axis, we get that the points $X, Y$ and $Z$ are collinear, belonging to the radical axis of the radical circle and to the circle circumscribed to the given triangle.

## 2nd Theorem.

If $M$ is a point in the plane of the triangle $ABC$ and the tangents in this point to the circles circumscribed to triangles $C, MAC, MAB$, respectively, cut $BC, CA$ and $AB$, respectively, in the points $X, Y, Z$, then these points are collinear.

## *Proof.*

The point $M$ is the radical center for the circles $(MBC)$, $(MAC)$, and $(MAB)$, and the tangents in $M$ to these circles are the polar lines to $M$ in relation to these circles.

If $X, Y, Z$ are the intersections of these tangents (polar lines) with $BC$, $CA, AB$, then they belong to the radical axis of the circumscribed circle to the triangle $ABC$ and to the circle "reduced" to the point $M$ ($XM^2 = XB \cdot XC$, etc.).

Being located on the radical axis of the two circles, the points $X, Y, Z$ are collinear.

## *Remarks.*

1. Another elementary proof of this theorem is to be found in [3].
2. If the circles from the $1^{st}$ *theorem* are adjoint circles of the triangle $ABC$ , then they intersect in $\Omega$ (the Brocard's point). Therefore, we get that the tangents taken in $\Omega$ to the adjoin circles cut the sides $BC$, $CA$ and $AB$ in collinear points.

# References.

# Regarding the First Droz-Farny's Circle

In this article, we define the **first Droz-Farny's circle**, we establish a connection between it and a **concyclicity theorem**, then we generalize this theorem, leading to the **generalization of Droz-Farny's circle**. The first theorem of this article was enunciated by J. Steiner and it was proven by Droz-Farny (*Mathésis*, 1901).

## 1$^{st}$ Theorem.

Let $ABC$ be a triangle, $H$ its orthocenter and $A_1, B_1, C_1$ the means of sides $(BC), (CA), (AB)$.

If a circle, having its center $H$, intersects $B_1C_1$ in $P_1, Q_1$ ; $C_1A_1$ in $P_2, Q_2$ and $A_1B_1$ in $P_3, Q_3$ , then $AP_1 = AQ_1 = BP_2 = BQ_2 = CP_3 = CQ_3$.

*Proof.*

Naturally, $HP_1 = HQ_1, B_1C_1 \parallel BC, AH \perp BC$.
It follows that $AH \perp B_1C_1$.

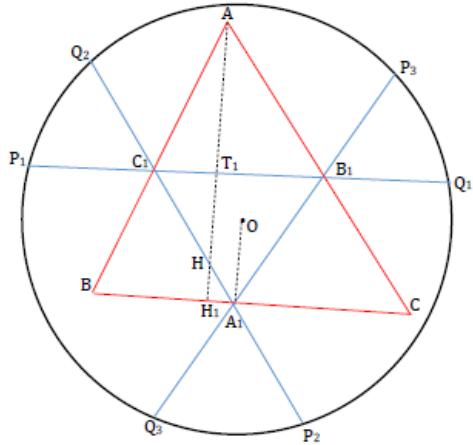

*Figure 1.*

Therefore, $AH$ is the mediator of segment $P_1Q_1$; similarly, $BH$ and $CH$ are the mediators of segments $P_2Q_2$ and $P_3Q_3$.

Let $T_1$ be the intersection of lines $AH$ and $B_1C_1$ (see *Figure 1*); we have $Q_1A^2 - Q_1H^2 = T_1A^2 - T_1H^2$. We denote $R_H = HP_1$. It follows that $Q_1A^2 = R_H^2 + (T_1A + T_1H)(T_1A - T_1H) = R_H^2 + AH \cdot (T_1A - T_1H)$.

However, $T_1A = T_1H_1$, where $H_1$ is the projection of $A$ on $BC$; we find that $Q_1A^2 = R_H^2 + AH \cdot HH_1$.

It is known that the symmetric of orthocenter $H$ towards $BC$ belongs to the circle of circumscribed triangle $ABC$.

Denoting this point by $H_1'$, we have $AH \cdot HH_1' = R^2 - OH^2$ (the power of point $H$ towards the circumscribed circle).

We obtain that $AH \cdot HH_1 = \frac{1}{2} \cdot (R^2 - OH^2)$, and therefore $AQ_1^2 = R_H^2 + \frac{1}{2} \cdot (R^2 - OH^2)$, where $O$ is the center of the circumscribed triangle $ABC$.

Similarly, we find $BQ_2^2 = CQ_3^2 = R_H^2 + \frac{1}{2}(R^2 - OH^2)$, therefore $AQ_1 = BQ_2 = CQ_3$.

*Remarks.*

- a. The proof adjusts analogously for the obtuse triangle.
- b. $1^{st}$ *Theorem* can be equivalently formulated in this way:

## $2^{nd}$ **Theorem.**

If we draw three congruent circles centered in a given triangle's vertices, the intersection points of these circles with the sides of the median triangle of given triangle (middle lines) are six points situated on a circle having its center in triangle's orthocenter.

If we denote by $\rho$ the radius of three congruent circles having $A, B, C$ as their centers, we get:

$R_H^2 = \rho^2 + \frac{1}{2}(OH^2 - R^2)$.

However, in a triangle, $OH^2 = 9R^2 - (a^2 + b^2 + c^2)$, $R$ being the radius of the circumscribed circle; it follows that:

$R_H^2 = \rho^2 + 4R^2 - \frac{1}{2}(a^2 + b^2 + c^2)$.

*Remark.*

A special case is the one in which $\rho = R$, where we find that $R_H^2 = R_1^2 = 5R^2 - \frac{1}{2}(a^2 + b^2 + c^2) = \frac{1}{2}(R^2 + OH^2)$.

## Definition.

The circle $\mathcal{C}(H, R_1)$, where:
$$R_1 = \sqrt{5R^2 - \frac{1}{2}(a^2 + b^2 + c^2)},$$
is called the first Droz-Farny's circle of the triangle $ABC$.

*Remark.*

Another way to build the first Droz-Farny's circle is offered by the following theorem, which, according to [1], was the subject of a problem proposed in 1910 by V. Thébault in the *Journal de Mathématiques Elementaire*.

## 3ʳᵈ Theorem.

The circles centered on the feet of a triangle's altitudes passing through the center of the circle circumscribed to the triangle cut the triangle's sides in six concyclical points.

### Proof.

We consider $ABC$ an acute triangle and $H_1, H_2, H_3$ the altitudes' feet. We denote by $A_1, A_2$; $B_1, B_2$; $C_1, C_2$ the intersection points of circles having their centers $H_1, H_2, H_3$ to $BC, CA, AB$, respectively.

We calculate $HA_2$ of the right angled triangle $HH_1A_2$ (see *Figure 2*). We have $HA_2^2 = HH_1^2 + H_1A_2^2$.

Because $H_1A_2 = H_1O$, it follows that $HA_2^2 = HH_1^2 + H_1O^2$. We denote by $O_9$ the mean of segment $OH$; the median theorem in triangle $H_1HO$ leads to $H_1H^2 + H_1O^2 = 2H_1O_9^2 + \frac{1}{2}OH^2$.

It is known that $O_9H_1$ is the nine-points circle's radius, so $H_1O_9 = \frac{1}{2}R$; we get: $HA_1^2 = \frac{1}{2}(R^2 + OH^2)$; similarly, we find that $HB_1^2 = HC_1^2 = \frac{1}{2}(R^2 + OH^2)$, which shows that the points $A_1, A_2$; $B_1, B_2$; $C_1, C_2$ belong to the first Droz-Farny's circle.

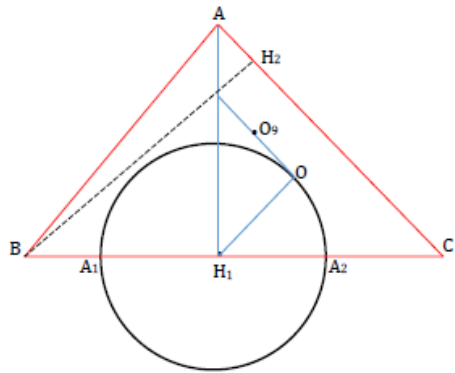

*Figure 2.*

*Remark.*

The $2^{nd}$ and the $3^{rd}$ theorems show that the first Droz-Farny's circle pass through 12 points, situated two on each side of the triangle and two on each side of the median triangle of the given triangle.

The following theorem generates the $3^{rd}$ *Theorem*.

# $4^{th}$ Theorem.

The circles centered on the feet of altitudes of a given triangle intersect the sides in six concyclic points if and only if their radical center is the center of the circle circumscribed to the given triangle.

*Proof.*

Let $A_1, A_2$; $B_1, B_2$; $C_1, C_2$ be the points of intersection with the sides of triangle $ABC$ of circles having their centers in altitudes' feet $H_1, H_2, H_3$.

Suppose the points are concyclic; it follows that their circle's radical center and of circles centered in $H_2$ and $H_3$ is the point $A$ (sides $AB$ and $AC$ are radical axes), therefore the perpendicular from $A$ on $H_2H_3$ is radical axis of centers having their centers $H_2$ and $H_3$.

Since $H_2H_3$ is antiparallel to $BC$, it is parallel to tangent taken in $A$ to the circle circumscribed to triangle $ABC$.

Consequently, the radical axis is the perpendicular taken in $A$ on the tangent to the circumscribed circle, therefore it is $AO$.

Similarly, the other radical axis of circles centered in $H_1, H_2$ and of circles centered in $H_1, H_3$ pass through $O$, therefore $O$ is the radical center of these circles.

### *Reciprocally.*

Let $O$ be the radical center of the circles having their centers in the feet of altitudes. Since $AO$ is perpendicular on $H_2H_3$, it follows that $AO$ is the radical axis of circles having their centers in $H_2, H_3$, therefore $AB_1 \cdot AB_2 = AC_1 \cdot AC_2$.

From this relationship, it follows that the points $B_1, B_2$; $C_1, C_2$ are concyclic; the circle on which these points are located has its center in the orthocenter $H$ of triangle $ABC$.

Indeed, the mediators' chords $B_1B_2$ and $C_1C_2$ in the two circles are the altitudes from $C$ and $B$ of triangle $ABC$, therefore $HB_1 = HB_2 = HC_1 = HC_2$.

This reasoning leads to the conclusion that $BO$ is the radical axis of circles having their centers $H_1$ and $H_3$, and from here the concyclicality of the points $A_1, A_2$; $C_1, C_2$ on a circle having its center in $H$, therefore $HA_1 = HA_2 = HC_1 = HC_2$. We obtained that $HA_1 = HA_2 = HB_1 = HB_2 = HC_1 = HC_2$, which shows the concyclicality of points $A_1, A_2, B_1, B_2, C_1, C_2$.

*Remark.*

The circles from the $3^{rd}$ *Theorem*, passing through *O* and having noncollinear centers, admit *O* as radical center, and therefore the $3^{rd}$ *Theorem* is a particular case of the $4^{th}$ *Theorem*.

# References.

# Regarding the Second Droz-Farny's Circle

In this article, we prove the theorem relative to the **second Droz-Farny's circle**, and a sentence that generalizes it.

The paper [1] informs that the following *Theorem* is attributed to J. Neuberg (*Mathesis*, 1911).

## 1$^{st}$ Theorem.

The circles with its centers in the middles of triangle $ABC$ passing through its orthocenter $H$ intersect the sides $BC$, $CA$ and $AB$ respectively in the points $A_1, A_2, B_1, B_2$ and $C_1, C_2$, situated on a concentric circle with the circle circumscribed to the triangle $ABC$ (the second Droz-Farny's circle).

### *Proof.*

We denote by $M_1, M_2, M_3$ the middles of $ABC$ triangle's sides, see *Figure 1*. Because $AH \perp M_2M_3$ and $H$ belongs to the circles with centers in $M_2$ and $M_3$, it follows that $AH$ is the radical axis of these circles, therefore we have $AC_1 \cdot AC_2 = AB_2 \cdot AB_1$. This relation shows that $B_1, B_2, C_1, C_2$ are concyclic points, because the center of the circle on which they are situated is $O$, the center of the circle circumscribed to the triangle $ABC$, hence we have that:

$$OB_1 = OC_1 = OC_2 = OB_2. \qquad (1)$$

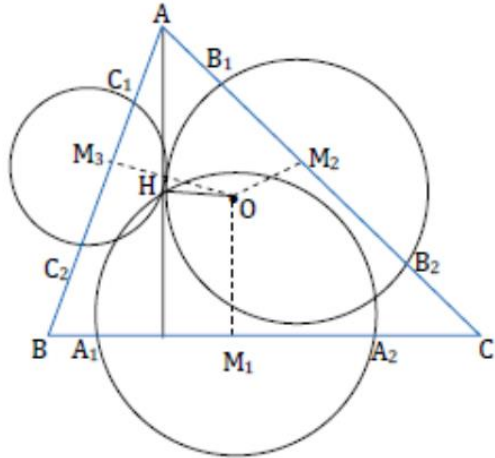

*Figure 1.*

Analogously, $O$ is the center of the circle on which the points $A_1, A_2, C_1, C_2$ are situated, hence:

$$OA_1 = OC_1 = OC_2 = OA_2. \qquad (2)$$

Also, $O$ is the center of the circle on which the points $A_1, A_2, B_1, B_2$ are situated, and therefore:

$$OA_1 = OB_1 = OB_2 = OA_2. \qquad (3)$$

The relations (1), (2), (3) show that the points $A_1, A_2, B_1, B_2, C_1, C_2$ are situated on a circle having the center in $O$, called the second Droz-Farny's circle.

# 1$^{\text{st}}$ Proposition.

The radius of the second Droz-Farny's circle is given by:

$$R_2^2 = 5e^2 - \frac{1}{2}(a^2 + b^2 + c^2).$$

## *Proof.*

From the right triangle $OM_1A_1$, using Pitagora's theorem, it follows that:

$$OA_1^2 = OM_1^2 + A_1M_1^2 = OM_1^2 + M_1M_2.$$

From the triangle $BHC$, using the median theorem, we have:

$$HM_1^2 = \frac{1}{4}[2(BH^2 + CH^2) - BC^2].$$

But in a triangle,

$$AH = 2OM_1,\ BH = 2OM_2,\ CH = 2OM_3,$$

hence:

$$HM_1^2 = 2OM_2^2 + 2OM_3^2 = \frac{a^2}{4}.$$

But:

$$OM_1^2 = R^2 - \frac{a^2}{4}\ ;$$
$$OM_2^2 = R^2 - \frac{b^2}{4}\ ;$$
$$OM_3^2 = R^2 - \frac{c^2}{4}\ ,$$

where R is the radius of the circle circumscribed to the triangle $ABC$.

We find that $OA_1^2 = R_2^2 = 5R^2 - \frac{1}{2}(a^2 + b^2 + c^2).$

*Remarks.*

a. We can compute $OM_1^2 + M_1M_2$ using the median theorem in the triangle $OM_1H$ for the median $M_1O_9$ ($O_9$ is the center of the nine points circle, i.e. the middle of $(OH)$). Because $O_9M_1 = \frac{1}{2}R$ , we obtain: $R_2^2 = \frac{1}{2}(OM^2 + R^2)$. In this way, we can prove the *Theorem* computing $OB_1^2$ and $OC_1^2$.

b. The statement of the *1st Theorem* was the subject no. 1 of the 49th International Olympiad in Mathematics, held at Madrid in 2008.

c. The *1st Theorem* can be proved in the same way for an obtuse triangle; it is obvious that for a right triangle, the second Droz-Farny's circle coincides with the circle circumscribed to the triangle $ABC$.

d. The *1st Theorem* appears as proposed problem in [2].

## 2nd Theorem.

The three pairs of points determined by the intersections of each circle with the center in the middle of triangle's side with the respective side are on a circle if and only these circles have as radical center the triangle's orthocenter.

*Proof.*

Let $M_1, M_2, M_3$ the middles of the sides of triangle $ABC$ and let $A_1, A_2, B_1, B_2, C_1, C_2$ the intersections with $BC$, $CA$, $AB$ respectively of the circles with centers in $M_1, M_2, M_3$.

Let us suppose that $A_1, A_2, B_1, B_2, C_1, C_2$ are concyclic points. The circle on which they are situated has evidently the center in $O$, the center of the circle circumscribed to the triangle $ABC$.

The radical axis of the circles with centers $M_2, M_3$ will be perpendicular on the line of centers $M_2 M_3$, and because A has equal powers in relation to these circles, since $AB_1 \cdot AB_2 = AC_1 \cdot AC_2$, it follows that the radical axis will be the perpendicular taken from A on $M_2 M_3$, i.e. the height from $A$ of triangle $ABC$.

Furthermore, it ensues that the radical axis of the circles with centers in $M_1$ and $M_2$ is the height from $B$ of triangle $ABC$ and consequently the intersection of the heights, hence the orthocenter $H$ of the triangle $ABC$ is the radical center of the three circles.

*Reciprocally.*

If the circles having the centers in $M_1, M_2, M_3$ have the orthocenter with the radical center, it follows that the point $A$, being situated on the height from A which is the radical axis of the circles of centers $M_2, M_3$ will have equal powers in relation to these circles and, consequently, $AB_1 \cdot AB_2 = AC_1 \cdot AC_2$ , a relation that implies that $B_1, B_2, C_1, C_2$ are concyclic points, and the circle on which these points are situated has $O$ as its center.

Similarly, $BA_1 \cdot BA_2 = BC_1 \cdot BC_2$ , therefore $A_1, A_2, C_1, C_2$ are concyclic points on a circle of center $O$. Having $OB_1 = OB_2 = OC_1 = OC_2$ and $OA_1 \cdot OA_2 = OC_1 \cdot OC_2$ , we get that the points $A_1, A_2, B_1, B_2, C_1, C_2$ are situated on a circle of center $O$.

*Remarks.*

1.  The *1<sup>st</sup> Theorem* is a particular case of the *2<sup>nd</sup> Theorem*, because the three circles of centers $M_1, M_2, M_3$ pass through $H$, which means that $H$ is their radical center.

2.  The Problem 525 from [3] leads us to the following *Proposition* providing a way to construct the circles of centers $M_1, M_2, M_3$ intersecting the sides in points that belong to a Droz-Farny's circle of type 2.

# 2<sup>nd</sup> Proposition.

The circles $\mathrm{C}\left(M_1, \frac{1}{2}\sqrt{k + a^2}\right)$, $\mathrm{C}\left(M_2, \frac{1}{2}\sqrt{k + b^2}\right)$, $\mathrm{C}\left(M_3, \frac{1}{2}\sqrt{k + c^2}\right)$ intersect the sides $BC$ , $CA$ , $AB$ respectively in six concyclic points; $k$ is a conveniently chosen constant, and $a, b, c$ are the lengths of the sides of triangle $ABC$.

*Proof.*

According to the *2$^{nd}$ Theorem*, it is necessary to prove that the orthocenter $H$ of triangle $ABC$ is the radical center for the circles from hypothesis.

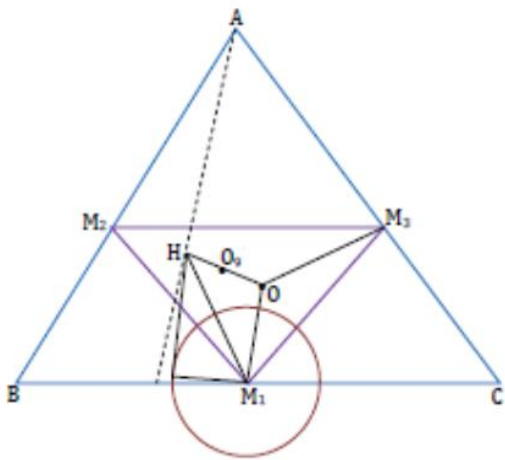

*Figure 2.*

The power of $H$ in relation with $C\left(M_1, \frac{1}{2}\sqrt{k+a^2}\right)$ is equal to $HM_1^2 - \frac{1}{4}(k+a^2)$ . We observed that $M_1^2 = 4R^2 - \frac{b^2}{2} - \frac{c^2}{2} - \frac{a^2}{4}$ , therefore $HM_1^2 - \frac{1}{4}(k+a^2) = 4R^2 - \frac{a^2+b^2+c^2}{4} - \frac{1}{4}k$. We use the same expression for the power of H in relation to the circles of centers $M_2, M_3$, hence H is the radical center of these three circles.

# References.

# Neuberg's Orthogonal Circles

In this article, we highlight some metric properties in connection with **Neuberg's circles** and **triangle**.

We recall some results that are necessary.

## 1st Definition.

It's called Brocard's point of the triangle $ABC$ the point $\Omega$ with the property: $\sphericalangle\Omega AB = \sphericalangle\Omega BC = \sphericalangle\Omega CA$. The measure of the angle $\Omega AB$ is denoted by $\omega$ and it is called Brocard's angle. It occurs the relationship:

$\operatorname{ctg}\omega = \operatorname{ctg}A + \operatorname{ctg}B + \operatorname{ctg}C$ (see [1]).

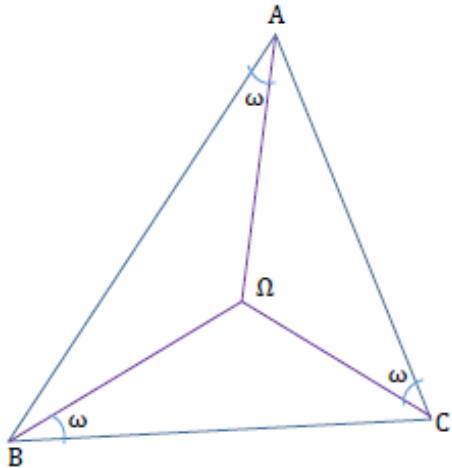

*Figure 1.*

## 2nd Definition.

Two triangles are called equibrocardian if they have the same Brocard's angle.

## 3rd Definition.

The locus of points $M$ from the plane of the triangle located on the same side of the line $BC$ as $A$ and forming with $BC$ an equibrocardian triangle with $ABC$, containing the vertex $A$ of the triangle, it's called A-Neuberg' circle of the triangle $ABC$.

We denote by $N_a$ the center of A-Neuberg' circle by radius $n_a$ (analogously, we define B-Neuberg' and C-Neuberg' circles).

We get that $m(BN_aC) = 2\omega$ and $n_a = \frac{a}{2}\sqrt{\text{ctg}^2\omega - 3}$ (see [1]).

The triangle $N_aN_bN_c$ formed by the centers of Neuberg's circles is called Neuberg's triangle.

## 1st Proposition.

The distances from the center circumscribed to the triangle $ABC$ to the vertex of Neuberg's triangle are proportional to the cubes of $ABC$ triangle's sides lengths.

*Proof.*

Let $O$ be the center of the circle circumscribed to the triangle $ABC$ (see *Figure 2*).

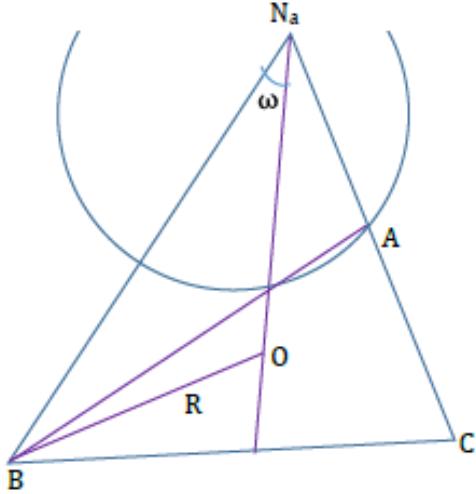

*Figure 2.*

The law of cosines applied in the triangle $ON_aB$ provides:

$\frac{ON_a}{\sin(N_aBO)} = \frac{R}{\sin\omega}$ .

But $m(\sphericalangle N_aBO) = m(\sphericalangle N_aBC) - m(\sphericalangle OBC) = A - \omega$.

We have that $\frac{ON_a}{\sin(A-\omega)} = \frac{R}{\sin\omega}$ .

But

$\frac{\sin(A-\omega)}{\sin\omega} = \frac{C\Omega}{A\Omega} = \frac{\frac{a}{c}\cdot 2R\sin\omega}{\frac{b}{a}\cdot 2R\sin\omega} = \frac{a^3}{abc} = \frac{a^3}{4RS}$ ,

$S$ being the area of the triangle $ABC$.

It follows that $ON_a = \frac{a^3}{4S}$, and we get that $\frac{ON_a}{a^3} = \frac{ON_b}{b^3} = \frac{ON_c}{c^3}$.

### *Consequence.*

In a triangle ABC, we have:
1)    $ON_a \cdot ON_b \cdot ON_c = R^3$;
2)    $\text{ctg}\omega = \frac{ON_a}{a} + \frac{ON_b}{b} + \frac{ON_c}{c}$.

## 2$^{\text{nd}}$ Proposition.

If $N_a N_b N_c$ is the Neuberg's triangle of the triangle $ABC$, then:

$$N_a N_b{}^2 = \frac{(a^2 + b^2)(a^4 + b^4) - a^2 b^2 c^2}{2a^2 b^2 + 2b^2 c^2 + 2c^2 a^2 - a^4 - b^4 - c^4}.$$

(The formulas for $N_b N_c$ and $N_c N_a$ are obtained from the previous one, by circular permutations.)

### *Proof.*

We apply the law of cosines in the triangle $N_a O N_c$:

$$ON_a = \frac{a^3}{4S} \, , ON_b = \frac{b^3}{4S} \, , m(\sphericalangle N_a O N_b) = 180^0 - \hat{c}.$$

$$N_a N_b{}^2 = \frac{a^6 + b^6 - 2a^3 b^3 \cos(180^0 - c)}{16S^2}$$

$$= \frac{a^6 + b^6 + 2a^3 b^3 \cos C}{16S^2}.$$

But the law of cosines in the triangle $ABC$ provides

$$2\cos C = \frac{a^2 + b^2 - c^2}{2ab}$$

and, from din Heron's formula, we find that

$$16S^2 = 2a^2b^2 + 2b^2c^2 + 2c^2a^2 - a^4 - b^4 - c^4.$$

Substituting the results above, we obtain, after a few calculations, the stated formula.

# 4th Definition.

Two circles are called orthogonal if they are secant and their tangents in the common points are perpendicular.

# 3rd Proposition.

(Gaultier – 18B)

Two circles $\mathcal{C}(O_1, r_1)$, $\mathcal{C}(O_2, r_2)$ are orthogonal if and only if

$$r_1^2 + r_2^2 = O_1O_2^2.$$

*Proof.*

Let $\mathcal{C}(O_1, r_1)$, $\mathcal{C}(O_2, r_2)$ be orthogonal (see *Figure 3*); then, if $A$ is one of the common points, the triangle $O_1AO_2$ is a right triangle and the Pythagorean Theorem applied to it, leads to $r_1^2 + r_2^2 = O_1O_2^2$.

### *Reciprocally.*

If the metric relationship from the statement occurs, it means that the triangle $O_1AO_2$ is a right triangle, therefore $A$ is their common point (the relationship $r_1{}^2 + r_2{}^2 = O_1O_2{}^2$ implies $r_1{}^2 + r_2{}^2 > O_1O_2{}^2$), then $O_1A \perp O_2A$, so $O_1A$ is tangent to the circle $\mathcal{C}(O_2, r_2)$ because it is perpendicular in $A$ on radius $O_2A$, and as well $O_2A$ is tangent to the circle $\mathcal{C}(O_1, r_1)$, therefore the circles are orthogonal.

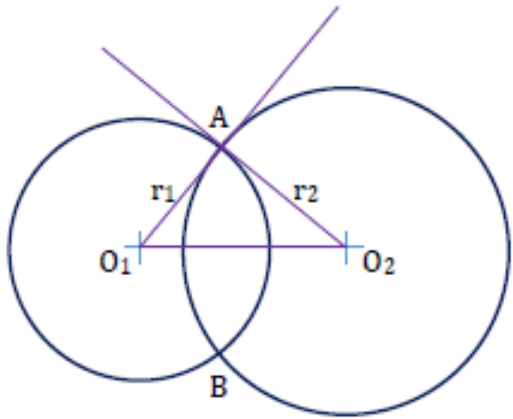

*Figure 3.*

# 4th Proposition.

B-Neuberg's and C-Neuberg's circles associated to the right triangle $ABC$ (in $A$) are orthogonal.

*Proof.*

If $m(\hat{A}) = 90^0$, then $N_b N_c^2 = \frac{b^6 + c^6}{16S}$.

$n_b = \frac{b}{2}\sqrt{\operatorname{ctg}^2\omega - 3}$ ; $n_c = \frac{c}{2}\sqrt{\operatorname{ctg}^2\omega - 3}$.

But $\operatorname{ctg}\omega - \frac{a^2 + b^2 + c^2}{4S} = \frac{b^2 + c^2}{2S} = \frac{a^2}{bc}$ .

It was taken into account that $a^2 = b^2 + c^2$ and $2S = bc$.

$$\operatorname{ctg}^2\omega - 3 = \frac{a^4}{b^2 c^2} - 3 = \frac{(b^2 + c^2)^2 - 3b^2 c^2}{b^2 c^2}$$

$$\operatorname{ctg}^2\omega - 3 = \frac{b^4 + c^4 - b^2 c^2}{b^2 c^2}$$

$$n_b^2 + n_c^2 = \frac{b^4 + c^4 - b^2 c^2}{b^2 c^2}\left(\frac{b^2 + c^2}{4}\right)$$

$$= \frac{(b^2 + c^2)(b^4 + c^4 - b^2 c^2)}{4b^2 c^2} = \frac{b^6 + c^6}{16S^2} .$$

By $N_b^2 + N_c^2 = N_b N_c^2$, it follows that B-Neuberg's and C-Neuberg's circles are orthogonal.

# References.

# Lucas's Inner Circles

In this article, we define the **Lucas's inner circles** and we highlight some of their properties.

## 1. Definition of the Lucas's Inner Circles

Let $ABC$ be a random triangle; we aim to construct the square inscribed in the triangle $ABC$, having one side on $BC$.

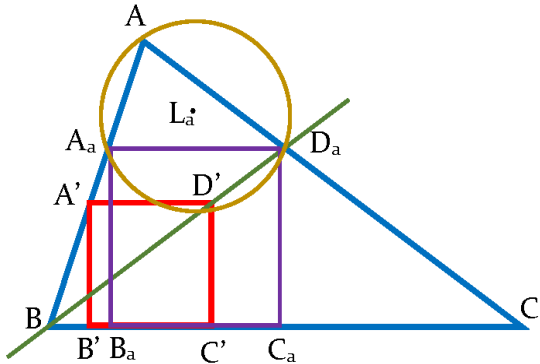

*Figure 1.*

In order to do this, we construct a square $A'B'C'D'$ with $A' \in (AB)$, $B', C' \in (BC)$ (see *Figure 1*).

We trace the line $BD'$ and we note with $D_a$ its intersection with $(AC)$; through $D_a$ we trace the parallel $D_aA_a$ to $BC$ with $A_a \in (AB)$ and we project onto $BC$ the points $A_a$, $D_a$ in $B_a$ respectively $C_a$.

We affirm that the quadrilateral $A_aB_aC_aD_a$ is the required square.

Indeed, $A_aB_aC_aD_a$ is a square, because $\frac{D_aC_a}{D'C'} = \frac{BD_a}{BD'} = \frac{A_aD_a}{A'D'}$ and, as $D'C' = A'D'$, it follows that $A_aD_a = D_aC_a$.

## Definition.

It is called A-Lucas's inner circle of the triangle $ABC$ the circle circumscribed to the triangle $AAaDa$.

We will note with $L_a$ the center of the A-Lucas's inner circle and with $l_a$ its radius.

Analogously, we define the B-Lucas's inner circle and the C-Lucas's inner circle of the triangle $ABC$.

## 2. Calculation of the Radius of the A-Lucas Inner Circle

We note $A_aD_a = x$, $BC = a$; let $h_a$ be the height from $A$ of the triangle $ABC$.

The similarity of the triangles $AA_aD_a$ and $ABC$ leads to: $\frac{x}{a} = \frac{h_a^{-x}}{h_a}$, therefore $x = \frac{ah_a}{a+h_a}$.

From $\frac{l_a}{R} = \frac{x}{a}$ we obtain $l_a = \frac{R.h_a}{a+h_a}$ . $\qquad$ (1)

*Note.*

Relation (1) and the analogues have been deduced by Eduard Lucas (1842-1891) in 1879 and they constitute the "birth certificate of the Lucas's circles".

$1^{st}$ *Remark.*

If in (1) we replace $h_a = \frac{2S}{a}$ and we also keep into consideration the formula $abc = 4RS$, where $R$ is the radius of the circumscribed circle of the triangle $ABC$ and $S$ represents its area, we obtain:

$l_a = \frac{R}{1+\frac{2aR}{bc}}$ [see Ref. 2].

## 3. Properties of the Lucas's Inner Circles

## $1^{st}$ Theorem.

The Lucas's inner circles of a triangle are inner tangents of the circle circumscribed to the triangle and they are exteriorly tangent pairwise.

*Proof.*

The triangles $AA_aD_a$ and $ABC$ are homothetic through the homothetic center $A$ and the rapport: $\frac{h_a}{a+h_a}$.

Because $\frac{l_a}{R} = \frac{h_a}{a + h_a}$, it means that the A-Lucas's inner circle and the circle circumscribed to the triangle $ABC$ are inner tangents in $A$.

Analogously, it follows that the B-Lucas's and C-Lucas's inner circles are inner tangents of the circle circumscribed to $ABC$.

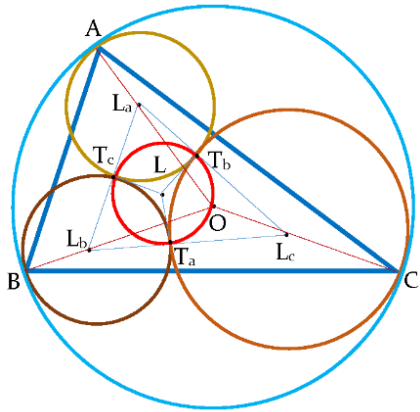

*Figure 2.*

We will prove that the A-Lucas's and C-Lucas's circles are exterior tangents by verifying

$$L_a L_c = l_a + l_c. \tag{2}$$

We have:

$$OL_a = R - l_a;$$
$$OL_c = R - l_c$$

and

$$m\left(\widehat{AOC}\right) = 2B$$

(if $m(\hat{B}) > 90°$ then $m\left(\widehat{AOC}\right) = 360° - 2B$).

The theorem of the cosine applied to the triangle $OL_aL_c$ implies, keeping into consideration (2), that:

$$(R - l_a)^2 + (R - l_a)^2 - 2(R - l_a)(R - l_c)cos2B =$$
$$= (l_a + l_c)^2.$$

Because $cos2B = 1 - 2sin^2B$, it is found that (2) is equivalent to:

$$sin^2B = \frac{l_al_c}{(R-l_a)(R-l_c)}. \tag{3}$$

But we have: $l_al_c = \frac{R^2ab^2c}{(2aR+bc)(2cR+ab)}$,

$l_a + l_c = Rb(\frac{c}{2aR+bc} + \frac{a}{2cR+ab})$.

By replacing in (3), we find that $sin^2B = \frac{ab^2c}{4acR^2} = \frac{b^2}{4a^2} \Leftrightarrow \sin B = \frac{b}{2R}$ is true according to the sines theorem. So, the exterior tangent of the A-Lucas's and C-Lucas's circles is proven.

Analogously, we prove the other tangents.

## 2$^{nd}$ Definition.

It is called an A-Apollonius's circle of the random triangle $ABC$ the circle constructed on the segment determined by the feet of the bisectors of angle $A$ as diameter.

### *Remark.*

Analogously, the B-Apollonius's and C-Apollonius's circles are defined. If $ABC$ is an isosceles triangle with $AB = AC$ then the A-Apollonius's circle isn't defined for $ABC$, and if $ABC$ is an equilateral triangle, its Apollonius's circle isn't defined.

## 2$^{nd}$ Theorem.

The A-Apollonius's circle of the random triangle is the geometrical point of the points $M$ from the plane of the triangle with the property: $\frac{MB}{MC} = \frac{c}{b}$.

## 3$^{rd}$ Definition.

We call a fascicle of circles the bunch of circles that do not have the same radical axis.

      a.      If the radical axis of the circles' fascicle is exterior to them, we say that the fascicle is of the first type.

      b.      If the radical axis of the circles' fascicle is secant to the circles, we say that the fascicle is of the second type.

      c.      If the radical axis of the circles' fascicle is tangent to the circles, we say that the fascicle is of the third type.

## 3$^{rd}$ Theorem.

The A-Apollonius's circle and the B-Lucas's and C-Lucas's inner circles of the random triangle $ABC$ form a fascicle of the third type.

*Proof.*

Let $\{O_A\} = L_b L_c \cap BC$ (see *Figure 3*).

Menelaus's theorem applied to the triangle $OBC$ implies that:

$$\frac{O_A B}{O_A C} \cdot \frac{L_b B}{L_b O} \cdot \frac{L_c O}{L_c C} = 1,$$

so:

$$\frac{O_A B}{O_A C} \cdot \frac{l_b}{R - l_b} \cdot \frac{R - l_c}{l_c} = 1$$

and by replacing $l_b$ and $l_c$, we find that:

$$\frac{O_A B}{O_A C} = \frac{b^2}{c^2}.$$

This relation shows that the point $O_A$ is the foot of the exterior symmedian from $A$ of the triangle $ABC$ (so the tangent in $A$ to the circumscribed circle), namely the center of the A-Apollonius's circle.

Let $N_1$ be the contact point of the B-Lucas's and C-Lucas's circles. The radical center of the B-Lucas's, C-Lucas's circles and the circle circumscribed to the triangle $ABC$ is the intersection $T_A$ of the tangents traced in $B$ and in $C$ to the circle circumscribed to the triangle $ABC$.

It follows that $BT_A = CT_A = N_1 T_A$, so $N_1$ belongs to the circle $\mathcal{C}_A$ that has the center in $T_A$ and orthogonally cuts the circle circumscribed in $B$ and $C$. The radical axis of the B-Lucas's and C-Lucas's circles is $T_A N_1$, and $O_A N_1$ is tangent in $N_1$ to the circle $\mathcal{C}_A$. Considering the power of the point $O_A$ in relation to $\mathcal{C}_A$, we have:

$$O_A N_1{}^2 = O_A B . O_A C.$$

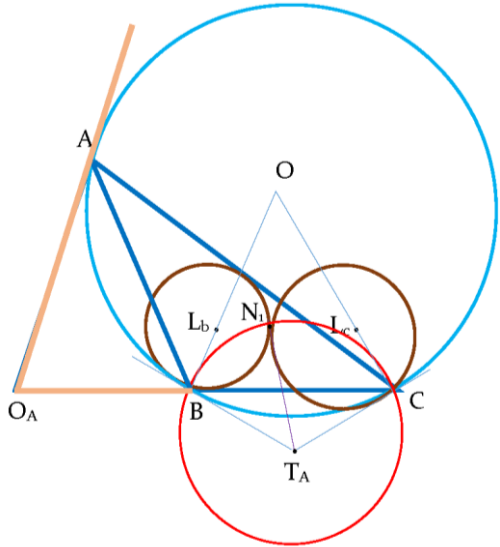

*Figure 3.*

Also, $O_A O^2 = O_A B \cdot O_A C$ ; it thus follows that $O_A A = O_A N_1$, which proves that $N_1$ belongs to the A-Apollonius's circle and is the radical center of the A-Apollonius's, B-Lucas's and C-Lucas's circles.

*Remarks.*

1.  If the triangle $ABC$ is right in $A$ then $L_b L_c || BC$, the radius of the A-Apollonius's circle is equal to: $\frac{abc}{|b^2 - c^2|}$. The point $N_1$ is the foot of the bisector from $A$. We find that $O_A N_1 = \frac{abc}{|b^2 - c^2|}$, so the theorem stands true.

2. The A-Apollonius's and A-Lucas's circles are orthogonal. Indeed, the radius of the A-Apollonius's circle is perpendicular to the radius of the circumscribed circle, $OA$, so, to the radius of the A-Lucas's circle also.

# 4$^{\text{th}}$ Definition.

The triangle $T_A T_B T_C$ determined by the tangents traced in $A, B, C$ to the circle circumscribed to the triangle $ABC$ is called the tangential triangle of the triangle $ABC$.

# 1$^{\text{st}}$ Property.

The triangle $ABC$ and the Lucas's triangle $L_a L_b L_c$ are homological.

*Proof.*

Obviously, $AL_a, BL_b, CL_c$ are concurrent in $O$, therefore $O$, the center of the circle circumscribed to the triangle $ABC$, is the homology center.

We have seen that $\{O_A\} = L_b L_c \cap BC$ and $O_A$ is the center of the A-Apollonius's circle, therefore the homology axis is the Apollonius's line $O_A O_B O_C$ (the line determined by the centers of the Apollonius's circle).

## 2$^{\text{nd}}$ **Property.**

The tangential triangle and the Lucas's triangle of the triangle $ABC$ are orthogonal triangles.

*Proof.*

The line $T_A N_1$ is the radical axis of the B-Lucas's inner circle and the C-Lucas's inner circle, therefore it is perpendicular on the line of the centers $L_b L_c$. Analogously, $T_B N_2$ is perpendicular on $L_c L_a$, because the radical axes of the Lucas's circles are concurrent in $L$, which is the radical center of the Lucas's circles; it follows that $T_A T_B T_C$ and $L_a L_b L_c$ are orthological and $L$ is the center of orthology. The other center of orthology is $O$ the center of the circle circumscribed to $ABC$.

## References.

# Theorems with Parallels Taken through a Triangle's Vertices and Constructions Performed only with the Ruler

In this article, we solve problems of **geometric constructions only with the ruler**, using known theorems.

## 1$^{st}$ Problem.

Being given a triangle $ABC$, its circumscribed circle (its center known) and a point $M$ fixed on the circle, construct, using only the ruler, a transversal line $A_1, B_1, C_1$, with $A_1 \in BC, B_1 \in CA, C_1 \in AB$, such that $\sphericalangle MA_1C \equiv \sphericalangle MB_1C \equiv \sphericalangle MC_1A$ (the lines taken though $M$ to generate congruent angles with the sides $BC$, $CA$ and $AB$, respectively).

## 2$^{nd}$ Problem.

Being given a triangle $ABC$, its circumscribed circle (its center known) and $A_1, B_1, C_1$, such that $A_1 \in$

$BC, B_1 \in CA, C_1 \in AB$ and $A_1, B_1, C_1$ collinear, construct, using only the ruler, a point $M$ on the circle circumscribing the triangle, such that the lines $MA_1, MB_1, MC_1$ to generate congruent angles with $BC$, $CA$ and $AB$, respectively.

## 3$^{\text{rd}}$ Problem.

Being given a triangle $ABC$ inscribed in a circle of given center and $AA'$ a given cevian, $A'$ a point on the circle, construct, using only the ruler, the isogonal cevian $AA_1$ to the cevian $AA'$.

To solve these problems and to prove the theorems for problems solving, we need the following *Lemma*:

## 1$^{\text{st}}$ Lemma.

(Generalized Simpson's Line)

If $M$ is a point on the circle circumscribed to the triangle $ABC$ and we take the lines $MA_1, MB_1, MC_1$ which generate congruent angles ( $A_1 \in BC, B_1 \in CA, C_1 \in AB$) with $BC, CA$ and $AB$ respectively, then the points $A_1, B_1, C_1$ are collinear.

*Proof.*

Let $M$ on the circle circumscribed to the triangle $ABC$ (see *Figure 1*), such that:

$$\sphericalangle MA_1C \equiv \sphericalangle MB_1C \equiv \sphericalangle MC_1A = \varphi. \qquad (1)$$

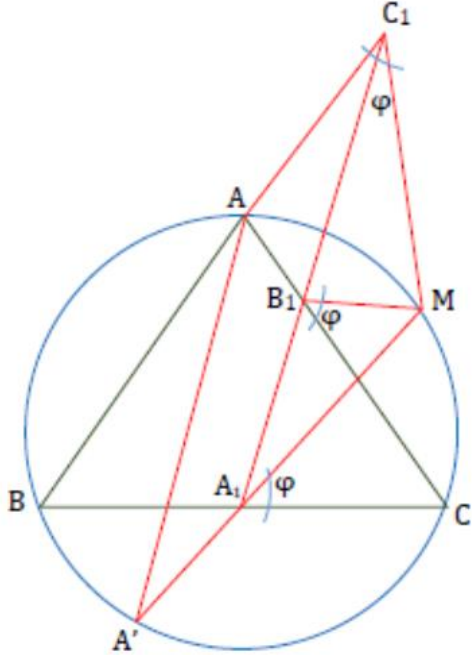

*Figure 1.*

From the relation (1), we obtain that the quadrilateral $MB_1A_1C$ is inscriptible and, therefore:

$$\sphericalangle A_1BC \equiv \sphericalangle A_1MC. \qquad (2).$$

Also from (1), we have that $MB_1AC_1$ is inscriptible, and so

$$\sphericalangle AB_1C_1 \equiv \sphericalangle AMC_1. \qquad (3)$$

The quadrilateral MABC is inscribed, hence:

$\sphericalangle MAC_1 \equiv \sphericalangle BCM.$       (4)

On the other hand,

$\sphericalangle A_1MC = 180^0 - \left(\widehat{BCM} + \varphi\right),$

$\sphericalangle AMC_1 = 180^0 - (\widehat{MAC_1} + \varphi).$

The relation (4) drives us, together with the above relations, to:

$\sphericalangle A_1MC \equiv \sphericalangle AMC_1.$       (5)

Finally, using the relations (5), (2) and (3), we conclude that: $\sphericalangle A_1B_1C \equiv AB_1C_1$, which justifies the collinearity of the points $A_1, B_1, C_1$.

### *Remark.*

The Simson's Line is obtained in the case when $\varphi = 90^0$.

## 2$^{\text{nd}}$ Lemma.

If $M$ is a point on the circle circumscribed to the triangle $ABC$ and $A_1, B_1, C_1$ are points on $BC$, $CA$ and $AB$, respectively, such that $\sphericalangle MA_1C = \sphericalangle MB_1C = \sphericalangle MC_1A = \varphi$, and $MA_1$ intersects the circle a second time in $A'$, then $AA' \parallel A_1B_1$.

### *Proof.*

The quadrilateral $MB_1A_1C$ is inscriptible (see *Figure 1*); it follows that:

$$\sphericalangle CMA' \equiv \sphericalangle A_1B_1C. \qquad (6)$$

On the other hand, the quadrilateral $MAA'C$ is also inscriptible, hence:

$$\sphericalangle CMA' \equiv \sphericalangle A'AC. \qquad (7)$$

The relations (6) and (7) imply: $\sphericalangle A'MC \equiv \sphericalangle A'AC$, which gives $AA' \parallel A_1B_1$.

## 3$^{\text{rd}}$ Lemma.

(The construction of a parallel with a given diameter using a ruler)

In a circle of given center, construct, using only the ruler, a parallel taken through a point of the circle at a given diameter.

### Solution.

In the given circle $\mathcal{C}(O, R)$, let be a diameter $(AB)]$ and let $M \in \mathcal{C}(O, R)$. We construct the line $BM$ (see *Figure 2*). We consider on this line the point $D$ ($M$ between $D$ and $B$). We join $D$ with $O$, $A$ with $M$ and denote $DO \cap AM = \{P\}$.

We take $BP$ and let $\{N\} = DA \cap BP$. The line $MN$ is parallel to $AB$.

### Construction's Proof.

In the triangle $DAB$, the cevians $DO$, $AM$ and $BN$ are concurrent.

Ceva's Theorem provides:

$$\frac{OA}{OB} \cdot \frac{MB}{MD} \cdot \frac{ND}{NA} = 1. \qquad (8)$$

But $DO$ is a median, $DO = BO = R$.

From (8), we get $\frac{MB}{MD} = \frac{NA}{ND}$ , which, by Thales reciprocal, gives $MN \parallel AB$.

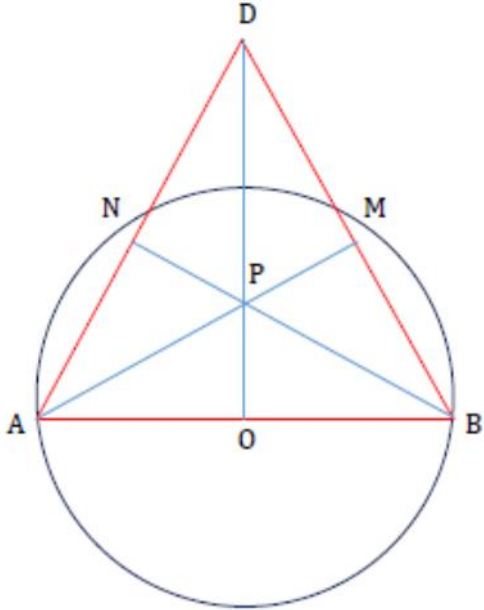

*Figure 2.*

### *Remark.*

If we have a circle with given center and a certain line $d$, we can construct though a given point $M$ a parallel to that line in such way: we take two diameters $[RS]$ and $[UV]$ through the center of the given circle (see *Figure 3*).

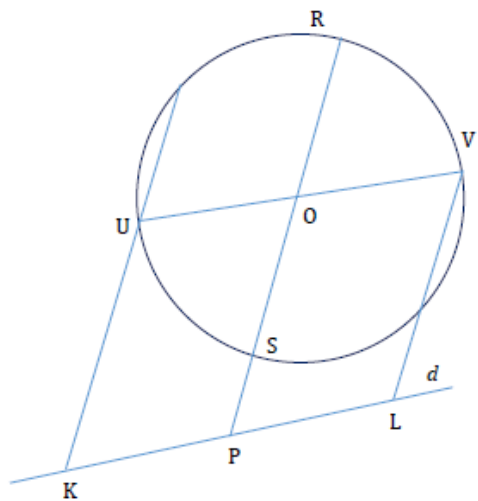

*Figure 3.*

We denote $RS \cap d = \{P\}$; because $[RO] \equiv [SO]$, we can construct, applying the $3^{rd}$ *Lemma*, the parallels through $U$ and $V$ to $RS$ which intersect $d$ in $K$ and $L$, respectively. Since we have on the line $d$ the points $K, P, L$, such that $[KP] \equiv [PL]$, we can construct the parallel through $M$ to $d$ based on the construction from $3^{rd}$ *Lemma*.

## 1st Theorem.

(P. Aubert – 1899)

If, through the vertices of the triangle $ABC$, we take three lines parallel to each other, which intersect the circumscribed circle in $A'$, $B'$ and $C'$, and $M$ is a point on the circumscribed circle, as well $MA' \cap BC = \{A_1\}$, $MB' \cap CA = \{B_1\}$, $MC' \cap AB = \{C_1\}$, then $A_1, B_1, C_1$ are collinear and their line is parallel to $AA'$.

## *Proof.*

The point of the proof is to show that $MA_1$, $MB_1$, $MC_1$ generate congruent angles with $BC$, $CA$ and $AB$, respectively.

$$m(\widehat{MA_1C}) = \frac{1}{2}\big[m(\widetilde{MC}) + m(\widetilde{BA'})\big] \qquad (9)$$

$$m(\widehat{MB_1C}) = \frac{1}{2}\big[m(\widetilde{MC}) + m(\widetilde{AB'})\big] \qquad (10)$$

But $AA' \parallel BB'$ implies $m(\widetilde{BA'}) = m(\widetilde{AB'})$, hence, from (9) and (10), it follows that:

$$\sphericalangle MA_1C \equiv \sphericalangle MB_1C, \qquad (11)$$

$$m(\widehat{MC_1A}) = \frac{1}{2}\big[m(\widetilde{BM}) - m(\widetilde{AC'})\big]. \qquad (12)$$

But $AA' \parallel CC'$ implies that $m(\widetilde{AC'}) = m(\widetilde{A'C})$; by returning to (12), we have that:

$$m(\widehat{MC_1A}) = \frac{1}{2}\big[m(\widetilde{BM}) - m(\widetilde{AC'})\big] =$$
$$= \frac{1}{2}\big[m(\widetilde{BA'}) + m(\widetilde{MC})\big]. \qquad (13)$$

The relations (9) and (13) show that:

$$\sphericalangle MA_1C \equiv \sphericalangle MC_1A. \qquad (14)$$

From (11) and (14), we obtain: $\sphericalangle MA_1C \equiv \sphericalangle MB_1C \equiv \sphericalangle MC_1A$, which, by *1st Lemma*, verifies the collinearity of points $A_1, B_1, C_1$. Now, applying the *2nd Lemma*, we obtain the parallelism of lines $AA'$ and $A_1B_1$.

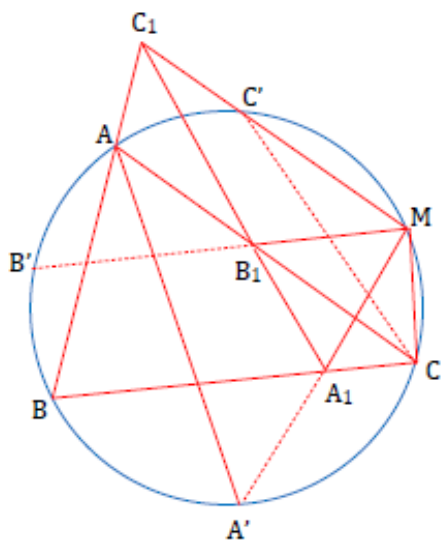

*Figure 4.*

## 2nd Theorem.

(M'Kensie – 1887)

If $A_1B_1C_1$ is a transversal line in the triangle $ABC$ ($A_1 \in BC, B_1 \in CA, C_1 \in AB$), and through the triangle's vertices we take the chords $AA'$, $BB'$, $CC'$ of a circle circumscribed to the triangle, parallels with the transversal line, then the lines $AA'$, $BB'$, $CC'$ are concurrent on the circumscribed circle.

### *Proof.*

We denote by $M$ the intersection of the line $A_1A'$ with the circumscribed circle (see *Figure 5*) and with $B_1'$, respectively $C_1'$ the intersection of the line $MB'$ with $AC$ and of the line $MC'$ with $AB$.

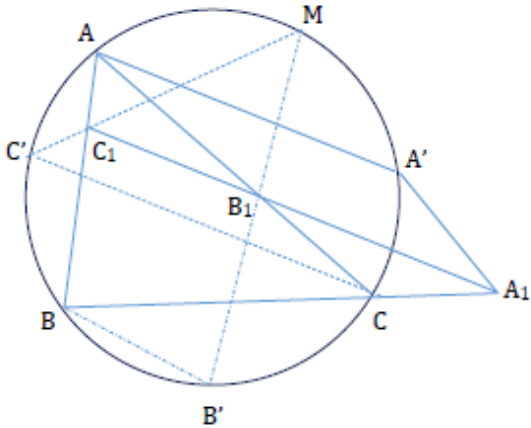

*Figure 5.*

According to the P. Aubert's theorem, we have that the points $A_1$, $B_1'$, $C_1'$ are collinear and that the line $A_1B_1'$ is parallel to $AA'$.

From hypothesis, we have that $A_1B_1 \parallel AA'$; from the uniqueness of the parallel taken through $A_1$ to $AA'$, it follows that $A_1B_1 \equiv A_1B_1'$, therefore $B_1' = B_1$, and analogously $C_1' = C_1$.

*Remark.*

We have that: $MA_1, MB_1, MC_1$ generate congruent angles with $BC$, $CA$ and $AB$, respectively.

## 3$^{\text{rd}}$ Theorem.

(Beltrami – 1862)

If three parallels are taken through the three vertices of a given triangle, then their isogonals intersect each other on the circle circumscribed to the triangle, and vice versa.

*Proof.*

Let $AA'$, $BB', CC'$ the three parallel lines with a certain direction (see *Figure 6*).

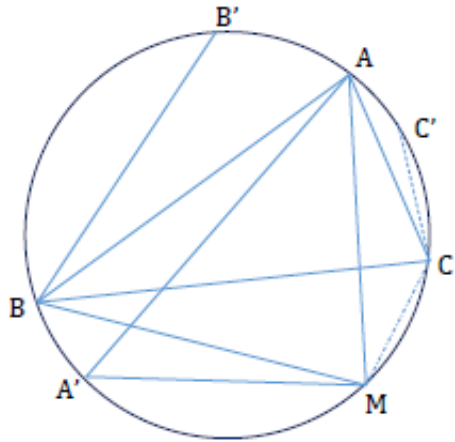

*Figure 6.*

To construct the isogonal of the cevian $AA'$, we take $A'M \parallel BC$, $M$ belonging to the circle circumscribed to the triangle, having $\widetilde{BA'} \equiv \widetilde{CM}$, it follows that $AM$ will be the isogonal of the cevian $AA'$. (Indeed, from $\widetilde{BA'} \equiv \widetilde{CM}$ it follows that $\sphericalangle BAA' \equiv \sphericalangle CAM$.)

On the other hand, $BB' \parallel$ A $A'$ implies $\widetilde{BA'} \equiv \widetilde{AB'}$, and since $\widetilde{BA'} \equiv \widetilde{CM}$ we have that $\widetilde{AB'} \equiv \widetilde{CM}$, which shows that the isogonal of the parallel $BB'$ is $BM$. From $CC' \parallel AA'$, it follows that $A'C \equiv AC'$, having $\sphericalangle B'CM \equiv \sphericalangle ACC'$, therefore the isogonal of the parallel $CC'$ is $CM'$.

### Reciprocally.

If $AM, BM, CM$ are concurrent cevians in $M$, the point on the circle circumscribed to the triangle $ABC$, let us prove that their isogonals are parallel lines. To construct an isogonal of $AM$, we take $MA' \parallel BC$, $A'$ belonging to the circumscribed circle. We have $\widetilde{MC} \equiv \widetilde{BA'}$. Constructing the isogonal $BB'$ of $BM$, with $B'$ on the circumscribed circle, we will have $\widetilde{CM} \equiv \widetilde{AB'}$, it follows that $\widetilde{BA'} \equiv \widetilde{AB'}$ and, consequently, $\sphericalangle ABB' \equiv \sphericalangle BAA'$, which shows that $AA' \parallel BB'$. Analogously, we show that $CC' \parallel AA'$.

We are now able to solve the proposed problems.

## Solution to the 1st problem.

Using the $3^{rd}$ *Lemma*, we construct the parallels $AA', BB', CC'$ with a certain directions of a diameter of the circle circumscribed to the given triangle.

We join $M$ with $A'$, $B', C'$ and denote the intersection between $MA'$ and $BC$, $A_1$; $MB' \cap CA = \{B_1\}$ and $MA' \cap AV = \{C_1\}$.

According to the Aubert's Theorem, the points $A_1, B_1, C_1$ will be collinear, and $MA'$, $MB'$, $MC'$ generate congruent angles with $BC$, $CA$ and $AB$, respectively.

## Solution to the 2nd problem.

Using the $3^{rd}$ *Lemma* and the remark that follows it, we construct through $A, B, C$ the parallels to $A_1B_1$; we denote by $A', B', C'$ their intersections with the circle circumscribed to the triangle $ABC$. (It is enough to build a single parallel to the transversal line $A_1B_1C_1$, for example $AA'$).

We join $A'$ with $A_1$ and denote by $M$ the intersection with the circle. The point $M$ will be the point we searched for. The construction's proof follows from the M'Kensie Theorem.

## Solution to the 3$^{rd}$ problem.

We suppose that $A'$ belongs to the little arc determined by the chord $\widetilde{BC}$ in the circle circumscribed to the triangle $ABC$.

In this case, in order to find the isogonal $AA_1$, we construct (by help of the $3^{rd}$ *Lemma* and of the remark that follows it) the parallel $A'A_1$ to $BC$, $A_1$ being on the circumscribed circle, it is obvious that $AA'$ and $AA_1$ will be isogonal cevians.

We suppose that $A'$ belongs to the high arc determined by the chord $\widetilde{BC}$; we consider $A' \in \widetilde{AB}$ (the arc $\widetilde{AB}$ does not contain the point $C$). In this situation, we firstly construct the parallel $BP$ to $AA'$, $P$ belongs to the circumscribed circle, and then through $P$ we construct the parallel $PA_1$ to $AC$, $A_1$ belongs to the circumscribed circle. The isogonal of the line $AA'$ will be $AA_1$. The construction's proof follows from $3^{rd}$ *Lemma* and from the proof of Beltrami's Theorem.

## References.

# Apollonius's Circles
# of $k^{\text{th}}$ Rank

The purpose of this article is to introduce the notion of **Apollonius's circle of $k^{\text{th}}$ rank**.

## 1<sup>st</sup> Definition.

It is called an internal cevian of $k^{th}$ rank the line $AA_k$ where $A_k \in (BC)$, such that $\frac{BA}{A_kC} = \left(\frac{AB}{AC}\right)^k$ $(k \in \mathbb{R})$.

If $A'_k$ is the harmonic conjugate of the point $A_k$ in relation to $B$ and $C$, we call the line $AA'_k$ an external cevian of $k^{th}$ rank.

## 2<sup>nd</sup> Definition.

We call Apollonius's circle of $k^{th}$ rank with respect to the side $BC$ of $ABC$ triangle the circle which has as diameter the segment line $A_k A'_k$.

## 1<sup>st</sup> Theorem.

Apollonius's circle of $k^{th}$ rank is the locus of points $M$ from $ABC$ triangle's plan, satisfying the relation: $\frac{MB}{MC} = \left(\frac{AB}{AC}\right)^k$.

## *Proof.*

Let $O_{A_k}$ the center of the Apollonius's circle of $k^{th}$ rank relative to the side $BC$ of $ABC$ triangle (see *Figure 1*) and $U, V$ the points of intersection of this circle with the circle circumscribed to the triangle $ABC$. We denote by $D$ the middle of arc $BC$, and we extend $DA_k$ to intersect the circle circumscribed in $U'$.

In $BU'C$ triangle, $U'D$ is bisector; it follows that $\frac{BA_k}{A_kC} = \frac{U'B}{U'C} = \left(\frac{AB}{AC}\right)^k$, so $U'$ belongs to the locus.

The perpendicular in $U'$ on $U'A_k$ intersects BC on $A_k''$, which is the foot of the $BUC$ triangle's outer bisector, so the harmonic conjugate of $A_k$ in relation to $B$ and $C$, thus $A_k'' = A_k'$.

Therefore, $U'$ is on the Apollonius's circle of rank $k$ relative to the side $BC$, hence $U' = U$.

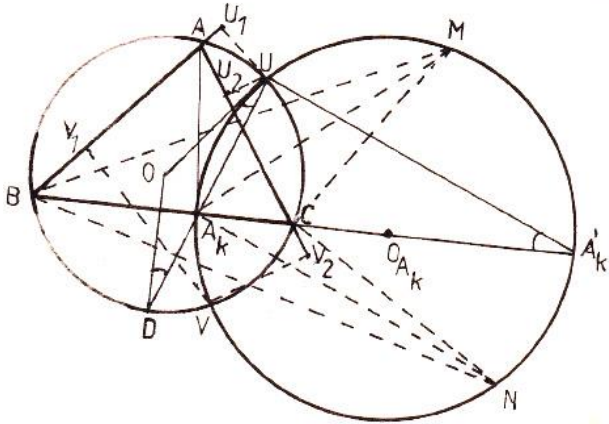

*Figure 3*

Let $M$ a point that satisfies the relation from the statement; thus $\frac{MB}{MC} = \frac{BA_k}{A_kC}$ ; it follows – by using the reciprocal of bisector's theorem – that $MA_k$ is the internal bisector of angle $BMC$. Now let us proceed as before, taking the external bisector; it follows that $M$ belongs to the Apollonius's circle of center $O_{A_k}$. We consider now a point $M$ on this circle, and we construct $C'$ such that $\sphericalangle BNA_k \equiv \sphericalangle A_kNC'$ (thus $(NA_k$ is the internal bisector of the angle $\widehat{BNC'}$ ). Because $A_k'N \perp NA_k$, it follows that $A_k$ and $A_k'$ are harmonically conjugated with respect to $B$ and $C'$. On the other hand, the same points are harmonically conjugated with respect to $B$ and $C$; from here, it follows that $C' = C$, and we have $\frac{NB}{NC} = \frac{BA_k}{A_kC} = \left(\frac{AB}{AC}\right)^k$.

# 3$^{\text{rd}}$ Definition.

It is called a complete quadrilateral the geometric figure obtained from a convex quadrilateral by extending the opposite sides until they intersect. A complete quadrilateral has 6 vertices, 4 sides and 3 diagonals.

# 2$^{\text{nd}}$ Theorem.

In a complete quadrilateral, the three diagonals' middles are collinear (Gauss - 1810).

### Proof.

Let $ABCDEF$ a given complete quadrilateral (see *Figure 2*). We denote by $H_1, H_2, H_3, H_4$ respectively the orthocenters of $ABF$, $ADE$, $CBE$, $CDF$ triangles, and let $A_1, B_1, F_1$ the feet of the heights of $ABF$ triangle.

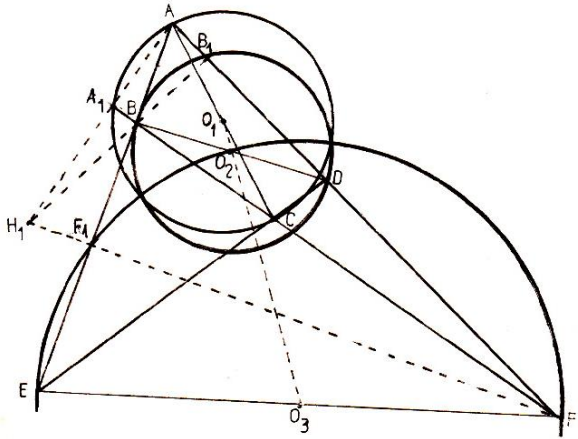

*Figure 4*

As previously shown, the following relations occur: $H_1A.H_1A_1 - H_1B.H_1B_1 = H_1F.H_1F_1$; they express that the point $H_1$ has equal powers to the circles of diameters $AC, BD, EF$, because those circles contain respectively the points $A_1, B_1, F_1$, and $H_1$ is an internal point.

It is shown analogously that the points $H_2, H_3, H_4$ have equal powers to the same circles, so those points are situated on the radical axis (common to the circles), therefore the circles are part of a fascicle, as such their centers – which are the middles of the complete quadrilateral's diagonals – are collinear.

The line formed by the middles of a complete quadrilateral's diagonals is called Gauss's line or Gauss-Newton's line.

## $3^{\text{rd}}$ Theorem.

The Apollonius's circle of $k^{th}$ rank of a triangle are part of a fascicle.

### *Proof.*

Let $AA_k$, $BB_k$, $CC_k$ be concurrent cevians of $k^{th}$ rank and $AA'_k$, $BB'_k$, $CC'_k$ be the external cevians of $k^{th}$ rank (see *Figure 3*). The figure $B'_k C_k B_k C'_k A_k A'_k$ is a complete quadrilateral and $2^{\text{nd}}$ theorem is applied.

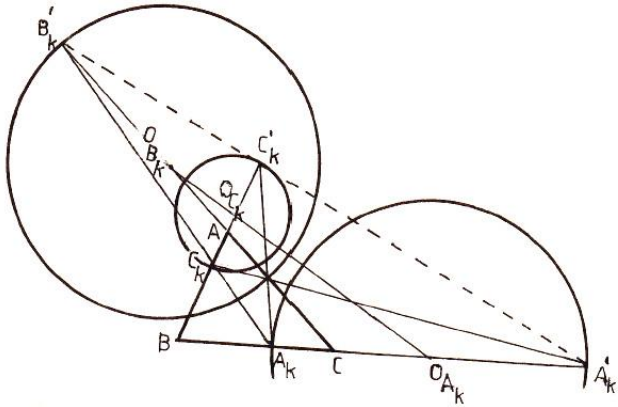

*Figure 5*

# 4$^{\text{th}}$ Theorem.

The Apollonius's circle of $k^{th}$ rank of a triangle are the orthogonals of the circle circumscribed to the triangle.

### *Proof.*

We unite $O$ to $D$ and $U$ (see *Figure 1*), $OD \perp BC$ and $m(\widehat{A_k U A'_k}) = 90^0$, it follows that $\widehat{U A'_k A_k} = \widehat{O D A_k} = \widehat{O U A_k}$.

The congruence $\widehat{U A'_k A_k} \equiv \widehat{O U A_k}$ shows that $OU$ is tangent to the Apollonius's circle of center $O_{A_k}$.

Analogously, it can be demonstrated for the other Apollonius's Circle.

### *1$^{st}$ Remark.*

The previous theorem indicates that the radical axis of Apollonius's circle of $k^{th}$ rank is the perpendicular taken from $O$ to the line $O_{A_k} O_{B_k}$.

# 5$^{\text{th}}$ Theorem.

The centers of Apollonius's Circle of $k^{th}$ rank of a triangle are situated on the trilinear polar associated to the intersection point of the cevians of $2k^{th}$ rank.

*Proof.*

From the previous theorem, it results that $OU \perp UO_{A_k}$, so $UO_{A_k}$ is an external cevian of rank 2 for $BCU$ triangle, thus an external symmedian. Henceforth, $\frac{O_{A_k}B}{O_{A_k}C} = \left(\frac{BU}{CU}\right)^2 = \left(\frac{AB}{AC}\right)^{2k}$ (the last equality occurs because $U$ belong to the Apollonius's circle of rank $k$ associated to the vertex $A$).

# $6^{\text{th}}$ Theorem.

The Apollonius's circle of $k^{th}$ rank of a triangle intersects the circle circumscribed to the triangle in two points that belong to the internal and external cevians of $k+1^{th}$ rank.

*Proof.*

Let $U$ and $V$ points of intersection of the Apollonius's circle of center $O_{A_k}$ with the circle circumscribed to the $ABC$ (see *Figure 1*). We take from $U$ and $V$ the perpendiculars $UU_1$, $UU_2$ and $VV_1, VV_2$ on $AB$ and $AC$ respectively. The quadrilaterals $ABVC$ , $ABCU$ are inscribed, it follows the similarity of triangles $BVV_1$, $CVV_2$ and $BUU_1$, $CUU_2$, from where we get the relations:

$$\frac{BV}{CV} = \frac{VV_1}{VV_2}, \qquad \frac{UB}{UC} = \frac{UU_1}{UU_2}.$$

But $\frac{BV}{CV} = \left(\frac{AB}{AC}\right)^k$, $\frac{UB}{UC} = \left(\frac{AB}{AC}\right)^k$, $\frac{VV_1}{VV_2} = \left(\frac{AB}{AC}\right)^k$ and $\frac{UU_1}{UU_2} = \left(\frac{AB}{AC}\right)^k$, relations that show that $V$ and $U$ belong respectively to the internal cevian and the external cevian of rank $k + 1$.

# 4$^{\text{th}}$ Definition.

If the Apollonius's circle of $k^{th}$ rank associated with a triangle has two common points, then we call these points isodynamic points of $k^{th}$ rank (and we denote them $W_k, W_k'$).

# 1$^{\text{st}}$ Property.

If $W_k, W_k'$ are isodynamic centers of $k^{th}$ rank, then:

$W_k A . BC^k = W_k B . AC^k = W_k C . AB^k$;
$W_k' A . BC^k = W_k' B . AC^k = W_k' C . AB^k$.

The proof of this property follows immediately from *1$^{st}$ Theorem*.

## *2$^{nd}$ Remark.*

The Apollonius's circle of 1$^{\text{st}}$ rank is the investigated Apollonius's circle (the bisectors are cevians of 1$^{\text{st}}$ rank). If $k = 2$, the internal cevians of 2$^{\text{nd}}$ rank are the symmedians, and the external cevians of 2$^{\text{nd}}$ rank are the external symmedians, i.e. the tangents in triangle's vertices to the circumscribed circle. In this case, for the Apollonius's circle of $2^{nd}$ rank, the *3rd Theorem* becomes:

# $7^{th}$ Theorem.

The Apollonius's circle of $2^{nd}$ rank intersects the circumscribed circle to the triangle in two points belonging respectively to the antibisector's isogonal and to the cevian outside of it.

*Proof.*

It follows from the proof of the $6^{th}$ theorem. We mention that the antibisector is isotomic to the bisector, and a cevian of $3^{rd}$ rank is isogonic to the antibisector.

# References.

# Apollonius's Circle of Second Rank

This article highlights some properties of **Apollonius's circle of second rank** in connection with the **adjoint circles** and the **second Brocard's triangle**.

## 1st Definition.

It is called Apollonius's circle of second rank relative to the vertex $A$ of the triangle $ABC$ the circle constructed on the segment determined on the simedians' feet from $A$ on $BC$ as diameter.

## 1st Theorem.

The Apollonius's circle of second rank relative to the vertex $A$ of the triangle $ABC$ intersect the circumscribed circle of the triangle $ABC$ in two points belonging respectively to the cevian of third rank (antibisector's isogonal) and to its external cevian.

The theorem's proof follows from the theorem relative to the Apollonius's circle of $k^{th}$ rank (see [1]).

# 1<sup>st</sup> Proposition.

The Apollonius's circle of second rank relative to the vertex $A$ of the triangle $ABC$ intersects the circumscribed circle in two points $Q$ and $P$ ($Q$ on the same side of $BC$ as $A$). Then, ($QS$ is a bisector in the triangle $QBC$, S is the simedian's foot from $A$ of the triangle $ABC$.

### *Proof.*

$Q$ belongs to the Apollonius's circle of second rank, therefore:

$$\frac{QB}{QC} = \left(\frac{AB}{AC}\right)^2. \qquad (1)$$

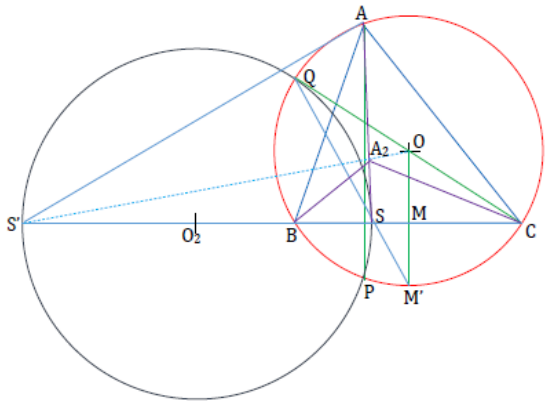

*Figure 1.*

On the other hand, $S$ being the simedian's foot, we have:

$\frac{SB}{SC} = \left(\frac{AB}{AC}\right)^2.$ \hfill (2)

From relations (1) and (2), we note that

$\frac{QB}{QC} = \frac{SB}{SC}$ ,

a relation showing that $QS$ is bisector in the triangle $QBC$.

### *Remarks.*

1.  The Apollonius's circle of second relative to the vertex $A$ of the triangle $ABC$ (see *Figure 1*) is an Apollonius's circle for the triangle $QBC$. Indeed, we proved that $QS$ is an internal bisector in the triangle $QBC$, and since $S'$, the external simedian's foot of the triangle $ABC$ , belongs to the Apollonius's Circle of second rank, we have $m(\sphericalangle S'QS) = 90^0$, therefore $QS'$ is an external bisector in the triangle $QBC$.

2.  $QP$ is a simedian in $QBC$ . Indeed, the Apollonius's circle of second rank, being an Apollonius's circle for $QBC$, intersects the circle circum-scribed to $QBC$ after $QP$, which is simedian in this triangle.

## 2$^{nd}$ Definition.

It is called adjoint circle of a triangle the circle that passes through two vertices of the triangle and in one of them is tangent to the triangle's side. We denote $(B\bar{A})$ the adjoint circle that passes through $B$ and $A$, and is tangent to the side $AC$ in $A$.

About the circles $(B\bar{A})$ and $(C\bar{A})$, we say that they are adjoint to the vertex $A$ of the triangle $ABC$.

## 3rd Definition.

It is called the second Brocard's triangle the triangle $A_2B_2C_2$ whose vertices are the projections of the center of the circle circumscribed to the triangle $ABC$ on triangle's simedians.

## 2nd Proposition.

The Apollonius's circle of second rank relative to the vertex $A$ of triangle $ABC$ and the adjoint circles relative to the same vertex $A$ intersect in vertex $A_2$ of the second Brocard's triangle.

*Proof.*

It is known that the adjoint circles $(B\bar{A})$ and $(C\bar{A})$ intersect in a point belonging to the simedian $AS$; we denote this point $A_2$ (see [3]).

We have:

$$\angle BA_2S = \angle A_2BA + \angle A_2AB,$$

but:

$$\angle A_2BA \equiv \angle BA_2S = \angle A_2AB + A_2AC = \angle A.$$

Analogously, $\sphericalangle CA_2 S = \sphericalangle A$, therefore $(A_2 S$ is the bisector of the angle $BA_2C$. The bisector's theorem in this triangle leads to:

$$\frac{SB}{SC} = \frac{BA_2}{CA_2} \, ,$$

but:

$$\frac{SB}{SC} = \left(\frac{AB}{AC}\right)^2 ,$$

consequently:

$$\frac{BA_2}{CA_2} = \left(\frac{AB}{AC}\right)^2 ,$$

so $A_2$ is a point that belongs to the Apollonius's circle of second rank.

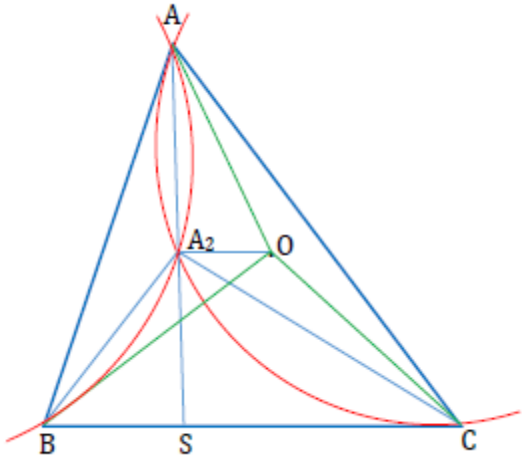

*Figure 2.*

We prove that $A_2$ is a vertex in the second Brocard's triangle, i.e. $OA_2 \perp AS$, O the center of the circle circumscribed to the triangle $ABC$.

We pointed (see *Figure 2*) that $m(\widehat{BA_2C}) = 2A$, if $\measuredangle A$ is an acute angle, then also $m(\widehat{BOC}) = 2A$, therefore the quadrilateral $OCA_2B$ is inscriptible.

Because $m(\widehat{OCB}) = 90^0 - m(\hat{A})$, it follows that $m(\widehat{BA_2O}) = 90^0 + m(\hat{A})$.

On the other hand, $m(\widehat{AA_2B}) = 180^0 - m(\hat{A})$, so $m(\widehat{BA_2O}) + m(\widehat{AA_2B}) = 270^0$ and, consequently, $OA_2 \perp AS$.

*Remarks.*

1.   If $m(\hat{A}) < 90^0$ , then four remarkable circles pass through $A_2$: the two circles adjoint to the vertex $A$ of the triangle $ABC$, the circle circumscribed to the triangle $BOC$ (where $O$ is the center of the circumscribed circle) and the Apollonius's circle of second rank corresponding to the vertex $A$.

2.   The vertex $A_2$ of the second Brocard's triangle is the middle of the chord of the circle circumscribed to the triangle $ABC$ containing the simedian $AS$.

3.   The points $O$, $A_2$ and $S$' (the foot of the external simedian to $ABC$) are collinear. Indeed, we proved that $OA_2 \perp AS$; on the other hand, we proved that ($A_2S$ is an internal bisector in the triangle $BA_2C$, and since $S'A_2 \perp AS$, the outlined collinearity follows from the uniqueness of the perpendicular in $A_2$ on $AS$.

## Open Problem.

The Apollonius's circle of second rank relative to the vertex $A$ of the triangle $ABC$ intersects the circle circumscribed to the triangle $ABC$ in two points $P$ and $Q$ ($P$ and $A$ apart of $BC$).

We denote by $X$ the second point of intersection between the line $AP$ and the Apollonius's circle of second rank.

What can we say about $X$?

Is $X$ a remarkable point in triangle's geometry?

# References.

# A Sufficient Condition for the Circle of the 6 Points to Become Euler's Circle

In this article, we prove the theorem relative to the **circle of the 6 points** and, requiring on this circle to have three other remarkable triangle's points, we obtain the **circle of 9 points** (the Euler's Circle).

## 1st Definition.

It is called cevian of a triangle the line that passes through the vertex of a triangle and an opposite side's point, other than the two left vertices.

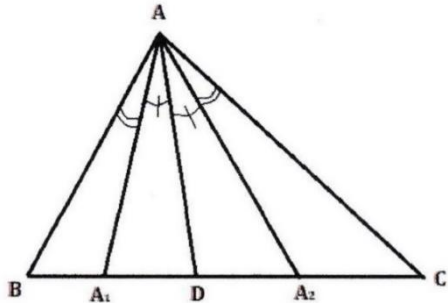

*Figure 1.*

## $1^{st}$ *Remark.*

The name *cevian* was given in honor of the italian geometrician Giovanni Ceva (1647 - 1734).

## $2^{nd}$ **Definition.**

Two cevians of the same triangle's vertex are called *isogonal cevians* if they create congruent angles with triangle's bisector taken through the same vertex.

## $2^{nd}$ *Remark.*

In the *Figure 1* we represented the bisector $AD$ and the isogonal cevians $AA_1$ and $AA_2$. The following relations take place:

$\widehat{A_1 AD} \equiv \widehat{A_2 AD}$;

$\widehat{BAA_1} \equiv \widehat{CAA_2}$.

## $1^{st}$ **Proposition.**

In a triangle, the height and the radius of the circumscribed circle corresponding to a vertex are isogonal cevians.

## *Proof.*

Let $ABC$ an acute-angled triangle with the height $AA'$ and the radius $AO$ (see *Figure 2*).

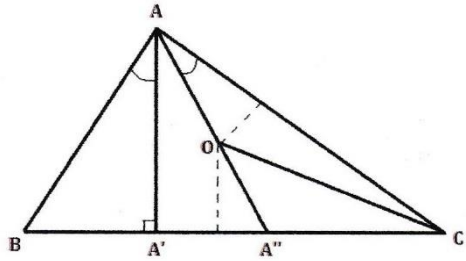

*Figure 2.*

The angle $AOC$ is a central angle, and the angle $ABC$ is an inscribed angle, so $\widehat{AOC} = 2\widehat{ABC}$. It follows that $\widehat{AOC} = 90^0 - \hat{B}$.

On the other hand, $\widehat{BAA'} = 90^0 - \hat{B}$, so $AA'$ and $AO$ are isogonal cevians.

The theorem can be analogously proved for the obtuse triangle.

### $3^{rd}$ *Remark.*

One can easily prove that in a triangle, if $AO$ is circumscribed circle's radius, its isogonal cevian is the triangle's height from vertex $A$.

## $3^{\text{rd}}$ Definition.

Two points $P_1, P_2$ in the plane of triangle $ABC$ are called *isogonal*s (isogonal conjugated) if the cevians' pairs $(AP_1, AP_2)$, $(BP_1, BP_2)$, $(CP_1, CP_2)$, are composed of isogonal cevians.

## 4<sup>th</sup> Remark.

In a triangle, the orthocenter and circumscribed circle's center are isogonal points.

## 1<sup>st</sup> Theorem.

(The 6 points circle)

If $P_1$ and $P_2$ are two isogonal points in the triangle $ABC$, and $A_1, B_1, C_1$ respectively $A_2, B_2, C_2$ are their projections on the sides $BC$, $CA$ and $AB$ of the triangle, then the points $A_1, A_2, B_1, B_2, C_1, C_2$ are concyclical.

### Proof.

The mediator of segment $[A_1 A_2]$ passes through the middle $P$ of segment $[P_1, P_2]$ because the trapezoid $P_1 A_1 A_2 P_2$ is rectangular and the mediator of $[A_1 A_2]$ contains its middle line, therefore (see *Figure 3*), we have: $PA_1 = PA_2$ (1). Analogously, it follows that the mediators of segments $[B_1 B_2]$ and $[C_1 C_2]$ pass through $P$, so $PB_1 = PB_2$ (2) and $PC_1 = PC_2$ (3). We denote by $A_3$ and $A_4$ respectively the middles of segments $[AP_1]$ and $[AP_2]$. We prove that the triangles $PA_3 C_1$ and $B_2 A_4 P$ are congruent. Indeed, $PA_3 = \frac{1}{2} AP_2$ (middle line), and $B_2 A_4 = \frac{1}{2} AP_2$, because it is a median in the rectangular triangle $P_2 B_2 A$, so $PA_3 = B_2 A_4$; analogously, we obtain that $A_4 P = A_3 C_1$.

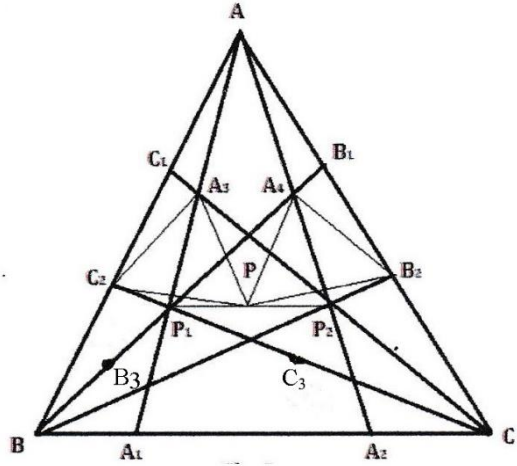

*Figure 3.*

We have that:

$$\widehat{PA_3C_1} = \widehat{PA_3P_1} + \widehat{P_1A_3C_1} = \widehat{P_1AP_2} + 2\widehat{P_1AC_1} =$$
$$= \hat{A} + \widehat{P_1AB};$$
$$\widehat{B_2A_4P} = \widehat{B_2A_4P_2} + \widehat{PA_4P_2} = \widehat{P_1AP_2} + 2\widehat{P_2AB_2} =$$
$$= \hat{A} + \widehat{P_2AC}.$$

But $\widehat{P_1AB} = \widehat{P_2AC}$, because the cevians $AP_1$ and $AP_2$ are isogonal and therefore $\widehat{PA_3C_1} = \widehat{B_2A_4P}$. Since $\Delta PA_3C_1 = \Delta B_2A_4P$, it follows that $PB_2 = PC_1$ (4).

Repeating the previous reasoning for triangles $PB_3C_1$ and $A_2B_4P$, where $B_3$ and $B_4$ are respectively the middles of segments $(BP_1)$ and $(BP_2)$, we find that they are congruent and it follows that $PC_1 = PA_2$ (5).

The relations (1) - (5) lead to $PA_1 = PA_2 = PB_1 = PB_2 = PC_1 = PC_2$ , which shows that the points $A_1, A_2, B_1, B_2, C_1, C_2$ are located on a circle centered in $P$, the middle of the segment determined by isogonal points given by $P_1$ and $P_2$.

# 4ᵗʰ Definition.

It is called the 9 points circle or Euler's circle of the triangle $ABC$ the circle that contains the middles of triangle's sides, the triangle heights' feet and the middles of the segments determined by the orthocenter with triangle's vertex.

# 2ⁿᵈ Proposition.

If $P_1, P_2$ are isogonal points in the triangle $ABC$ and if on the circle of their corresponding 6 points there also are the middles of the segments $(AP_1)$, $(BP_1)$, $(CP_1)$, then the 6 points circle coincides with the Euler's circle of the triangle $ABC$.

## *1ˢᵗ Proof.*

We keep notations from *Figure 3*; we proved that the points $A_1, A_2, B_1, B_2, C_1, C_2$ are on the 6 points circle of the triangle $ABC$, having its center in $P$, the middle of segment $[P_1P_2]$.

If on this circle are also situated the middles $A_3, B_3, C_3$ of segments $(AP_1)$, $(BP_1)$, $(CP_1)$, then we have $PA_3 = PB_3 = PC_3$.

We determined that $PA_3$ is middle line in the triangle $P_1P_2A$, therefore $PA_3 = \frac{1}{2}AP_2$, analogously $PB_3 = \frac{1}{2}BP_2$ and $PC_3 = \frac{1}{2}CP_2$, and we obtain that $P_2A = P_2B = P_2C$, consequently $P$ is the center of the circle circumscribed to the triangle $ABC$, so $P_2 = O$.

Because $P_1$ is the isogonal of $O$, it follows that $P_1 = H$, therefore the circle of 6 points of the isogonal points $O$ and $H$ is the circle of 9 points.

## $2^{nd}$ *Proof.*

Because $A_3B_3$ is middle line in the triangle $P_1AB$, it follows that

$$\sphericalangle P_1AB \equiv \sphericalangle P_1A_2B_3. \tag{1}$$

Also, $A_3C_3$ is middle line in the triangle $P_1AC$, and $A_3C_3$ is middle line in the triangle $P_1AP_2$, therefore we get

$$\sphericalangle PA_3C_3 \equiv \sphericalangle P_2AC. \tag{2}$$

The relations (1), (2) and the fact that $AP_1$ and $AP_2$ are isogonal cevians lead to:

$$\sphericalangle P_1A_2B_3 \equiv PA_3C_3. \tag{3}$$

The point $P$ is the center of the circle circumscribed to $A_3B_3C_3$; then, from (3) and from isogonal cevians' properties, one of which is circumscribed circle radius, it follows that in the triangle $A_3B_3C_3$ the line $P_1A_3$ is a height, as $B_3C_3 \parallel BC$, we get that $P_1A$ is a height in the triangle ABC and, therefore, $P_1$ will be the orthocenter of the triangle $ABC$ , and $P_2$ will be the center of the circle circumscribed to the triangle $ABC$.

## $5^{th}$ Remark.

Of those shown, we get that the center of the circle of 9 points is the middle of the line determined by triangle's orthocenter and by the circumscribed circle's center, and the fact that Euler's circle radius is half the radius of the circumscribed circle.

## References.

# An Extension of a Problem
# of Fixed Point

In this article, we extend the requirement of the *Problem 9.2* proposed at *Varna 2015 Spring Competition*, both in terms of membership of the measure $\gamma$, and the case for the problem for the ex-inscribed circle $C$. We also try to guide the student in the search and identification of the fixed point, for succeeding in solving any problem of this type.

The statement of the problem is as follows:

"We fix an angle $\gamma \in (0, 90^0)$ and the line $AB$ which divides the plane in two half-planes $\gamma$ and $\bar{\gamma}$. The point $C$ in the half-plane $\gamma$ is situated such that $m(\widehat{ACB}) = \gamma$. The circle inscribed in the triangle $ABC$ with the center $I$ is tangent to the sides $AC$ and $BC$ in the points $F$ and $E$, respectively. The point $P$ is located on the segment line $(IE$, the point $E$ between $I$ and $P$ such that $PE \perp BC$ and $PE = AF$. The point $Q$ is situated on the segment line $(IF$, such that $F$ is between $I$ and $Q$; $QF \perp AC$ and $QF = BE$. Prove that the mediator of segment $PQ$ passes through a fixed point." *(Stanislav Chobanov)*

*Proof.*

Firstly, it is useful to note that the point $C$ varies in the half-plane $\psi$ on the arc capable of angle $\gamma$; we know as well that $m\left(\widehat{AIB}\right) = 90^0 + \frac{\gamma}{2}$, so $I$ varies on the arc capable of angle of measure $90^0 + \frac{\gamma}{2}$ situated in the half-plane $\psi$.

Another useful remark is about the segments $AF$ and $BE$, which in a triangle have the lengths $p - a$, respectively $p - b$, where $p$ is the half-perimeter of the triangle $ABC$ with $AB = c$ – constant; therefore, we have $\Delta PEB \equiv \Delta AFQ$ with the consequence $PB = QA$. Considering the vertex $C$ of the triangle $ABC$ the middle of the arc capable of angle $\gamma$ built on $AB$, we observe that $PQ$ is parallel to $AB$; more than that, $ABPQ$ is an isosceles trapezoid, and segment $PQ$ mediator will be a symmetry axis of the trapezoid, so it will coincide with the mediator of $AB$, which is a fixed line, so we're looking for the fixed point on mediator of $AB$.

Let $D$ be the intersection of the mediators of segments $PQ$ and $AB$, see *Figure 1*, where we considered $m(\hat{A}) < m(\hat{B})$. The point $D$ is on the mediator of $AB$, so we have $DA = DB$; the point $D$ is also on the mediator of $PQ$, so we have $DP = DQ$; it follows that: $\Delta PBD \equiv \Delta QAD$, a relation from where we get that $\sphericalangle QAD = \sphericalangle PBD$.

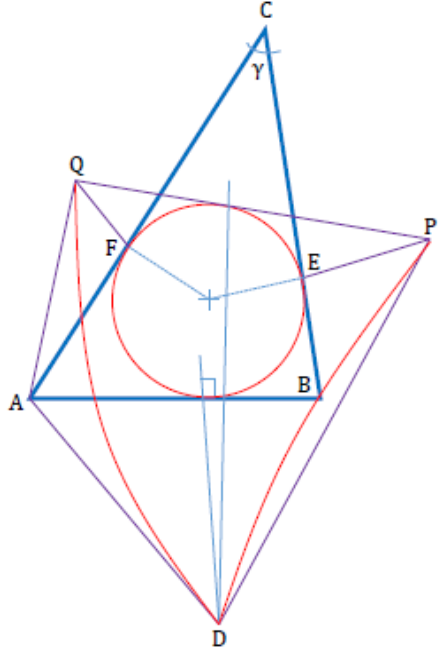

*Figure 1*

If we denote $m\left(\widehat{QAF}\right) = x$ and $m\left(\widehat{DAB}\right) = y$, we have $\widehat{QAD} = x + A + y$, $\widehat{PBD} = 360^0 - B - y - (90^0 - x)$.

From $x + A + y = 360^0 - B - y - 90^0 + x$, we find that $A + B + 2y = 270^0$, and since $A + B = 180^0 - \gamma$, we find that $2y = 90^0 - \gamma$, therefore the requested fixed point $D$ is the vertex of triangle $DAB$, situated in $\bar{\psi}$ such that $m\left(\widehat{ADB}\right) = 90^0 - \gamma$.

## $1^{st}$ **Remark.**

If $\gamma = 90^0$, we propose to the reader to prove that the quadrilateral $ABPQ$ is a parallelogram; in this case, the requested fixed point does not exist (or is the point at infinity of the perpendicular lines to $AB$).

## $2^{nd}$ **Remark.**

If $\gamma \in (90^0, 180^0)$, the problem applies, and we find that the fixed point $D$ is located in the half-plane $\psi$, such that the triangle $DAB$ is isosceles, having $m(\widehat{AOB}) = \gamma - 90^0$.

We suggest to the reader to solve the following problem:

> We fix an angle $\gamma \in (0^0, 180^0)$ and the line AB which divides the plane in two half-planes, $\psi$ and $\bar{\psi}$. The point $C$ in the half-plane $\psi$ is located such that $m(\widehat{ACB}) = \gamma$. The circle $C$ – ex-inscribed to the triangle $ABC$ with center $I_c$ is tangent to the sides $AC$ and $BC$ in the points $F$ and $E$, respectively. The point $P$ is located on the line segment $(I_c E$, $E$ is between $I_c$ and $P$ such that $PE \perp BC$ and $PE = AF$. The point $Q$ is located on the line segment $(I_c F$ such that $F$ is between $I$ and $Q$, $QF \perp AC$ and $QF = BE$. Prove that the mediator of the segment $PQ$ passes through a fixed point.

*3$^{rd}$ Remark.*

As seen, this problem is also true in the case $\gamma = 90^0$, more than that, in this case, the fixed point is the middle of $AB$. Prove!

## References.

# Some Properties of the Harmonic Quadrilateral

In this article, we review some properties of the **harmonic quadrilateral** related to triangle simedians and to Apollonius's Circle.

## 1<sup>st</sup> Definition.

A convex circumscribable quadrilateral $ABCD$ having the property $AB \cdot CD = BC \cdot AD$ is called harmonic quadrilateral.

## 2<sup>nd</sup> Definition.

A triangle simedian is the isogonal cevian of a triangle median.

## 1<sup>st</sup> Proposition.

In the triangle $ABC$, the cevian $AA_1$, $A_1 \in (BC)$ is a simedian if and only if $\frac{BA_1}{A_1C} = \left(\frac{AB}{AC}\right)^2$. For *Proof* of this property, see *infra*.

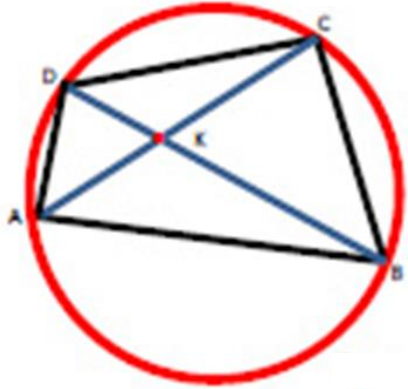

*Figura 1.*

## 2$^{\text{nd}}$ **Proposition.**

In an harmonic quadrilateral, the diagonals are simedians of the triangles determined by two consecutive sides of a quadrilateral with its diagonal.

### *Proof.*

Let $ABCD$ be an harmonic quadrilateral and $\{K\} = AC \cap BD$ (see *Figure 1*). We prove that $BK$ is simedian in the triangle $ABC$.

From the similarity of the triangles $ABK$ and $DCK$, we find that:

$$\frac{AB}{DC} = \frac{AK}{DK} = \frac{BK}{CK}. \tag{1}$$

From the similarity of the triangles $BCK$ şi $ADK$, we conclude that:

$$\frac{BC}{AD} = \frac{CK}{DK} = \frac{BK}{AK}. \tag{2}$$

From the relations (1) and (2), by division, it follows that:

$$\frac{AB}{BC} \cdot \frac{AD}{DC} = \frac{AK}{CK}. \tag{3}$$

But $ABCD$ is an harmonic quadrilateral; consequently,

$$\frac{AB}{BC} = \frac{AD}{DC} \; ;$$

substituting this relation in (3), it follows that:

$$\left(\frac{AB}{BC}\right)^2 = \frac{AK}{CK};$$

As shown by Proposition 1, $BK$ is a simedian in the triangle $ABC$. Similarly, it can be shown that $AK$ is a simedian in the triangle $ABD$, that $CK$ is a simedian in the triangle $BCD$, and that $DK$ is a simedian in the triangle $ADC$.

### *Remark 1.*

The converse of the 2$^{nd}$ Proposition is proved similarly, i.e.:

# 3$^{rd}$ **Proposition.**

If in a convex circumscribable quadrilateral, a diagonal is a simedian in the triangle formed by the other diagonal with two consecutive sides of the quadrilateral, then the quadrilateral is an harmonic quadrilateral.

*Remark 2.*

From 2<sup>nd</sup> and 3<sup>rd</sup> Propositions above, it results a simple way to build an harmonic quadrilateral.

In a circle, let a triangle $ABC$ be considered; we construct the simedian of $A$, be it $AK$, and we denote by $D$ the intersection of the simedian $AK$ with the circle. The quadrilateral $ABCD$ is an harmonic quadrilateral.

## Proposition 4.

In a triangle $ABC$, the points of the simedian of $A$ are situated at proportional lengths to the sides $AB$ and $AC$.

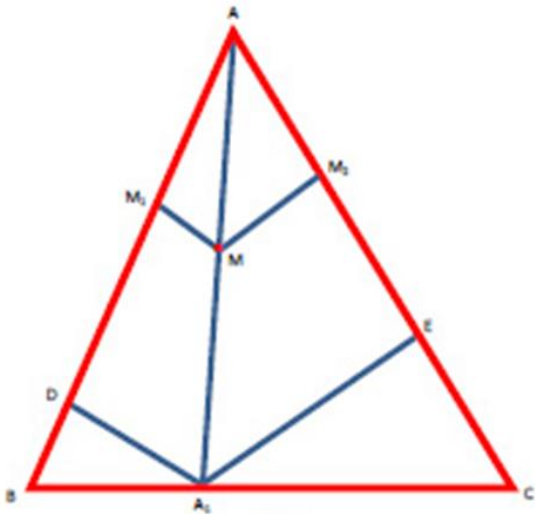

*Figura 2.*

*Proof.*

We have the simedian $AA_1$ in the triangle $ABC$ (see *Figure 2*). We denote by $D$ and $E$ the projections of $A_1$ on $AB$, and $AC$ respectively.

We get:

$$\frac{BA_1}{CA_1} = \frac{Area_\Delta(ABA_1)}{Area_\Delta(ACA_1)} = \frac{AB \cdot A_1D}{AC \cdot A_1E}.$$

Moreover, from $1^{st}$ Proposition, we know that

$$\frac{BA_1}{A_1C} = \left(\frac{AB}{AC}\right)^2.$$

Substituting in the previous relation, we obtain that:

$$\frac{A_1D}{A_1E} = \frac{AB}{AC}.$$

On the other hand, $DA_1 = AA_1$. From $BAA_1$ and $A_1E = AA_1 \cdot sin\widehat{CAA_1}$, hence:

$$\frac{A_1D}{A_1E} = \frac{sin\widehat{BAA_1}}{sin\widehat{CAA_1}} = \frac{AB}{AC}. \tag{4}$$

If $M$ is a point on the simedian and $MM_1$ and $MM_2$ are its projections on $AB$, and $AC$ respectively, we have:

$$MM_1 = AM \cdot sin\widehat{BAA_1}, \qquad MM_2 = AM \cdot sin\widehat{CAA_1},$$

hence:

$$\frac{MM_1}{MM_2} = \frac{sin\widehat{BAA_1}}{sin\widehat{CAA_1}}.$$

Taking (4) into account, we obtain that:

$$\frac{MM_1}{MM_2} = \frac{AB}{AC}.$$

## $3^{rd}$ *Remark.*

The converse of the property in the statement above is valid, meaning that, if  *M* is a point inside a triangle, its distances to two sides are proportional to the lengths of these sides. The point belongs to the simedian of the triangle having the vertex joint to the two sides.

## $5^{th}$ **Proposition.**

In an harmonic quadrilateral, the point of intersection of the diagonals is located towards the sides of the quadrilateral to proportional distances to the length of these sides. The *Proof* of this Proposition relies on $2^{nd}$ and $4^{th}$ Propositions.

## $6^{th}$ **Proposition.**

(R. Tucker)

The point of intersection of the diagonals of an harmonic quadrilateral minimizes the sum of squares of distances from a point inside the quadrilateral to the quadrilateral sides.

### *Proof.*

Let *ABCD* be an harmonic quadrilateral and *M* any point within.

We denote by $x, y, z, u$ the distances of $M$ to the $AB$ , $BC$ , $CD$ , $DA$ sides of lenghts $a, b, c,$ and $d$ (see *Figure 3*).

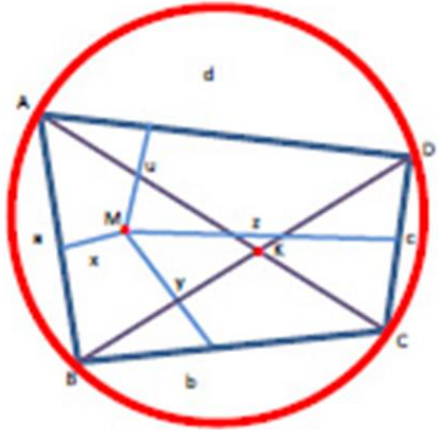

*Figure 3.*

Let $S$ be the $ABCD$ quadrilateral area.

We have:

$ax + by + cz + du = 2S.$

This is true for $x, y, z, u$ and $a, b, c, d$ real numbers.

Following Cauchy-Buniakowski-Schwarz Inequality, we get:

$$(a^2 + b^2 + c^2 + d^2)(x^2 + y^2 + z^2 + u^2)$$
$$\geq (ax + by + cz + du)^2 ,$$

and it is obvious that:

$$x^2 + y^2 + z^2 + u^2 \geq \frac{4S^2}{a^2 + b^2 + c^2 + d^2} .$$

We note that the minimum sum of squared distances is:

$$\frac{4S^2}{a^2 + b^2 + c^2 + d^2} = const.$$

In Cauchy-Buniakowski-Schwarz Inequality, the equality occurs if:

$$\frac{x}{a} = \frac{y}{b} = \frac{z}{c} = \frac{u}{d}.$$

Since $\{K\} = AC \cap BD$ is the only point with this property, it ensues that $M = K$, so $K$ has the property of the minimum in the statement.

## 3$^{\text{rd}}$ Definition.

We call external simedian of $ABC$ triangle a cevian $AA_1'$ corresponding to the vertex $A$, where $A_1'$ is the harmonic conjugate of the point $A_1$ – simedian's foot from $A$ relative to points $B$ and $C$.

### 4$^{th}$ Remark.

In *Figure 4,* the cevian $AA_1$ is an internal simedian, and $AA_1'$ is an external simedian.

We have:

$$\frac{A_1B}{A_1C} = \frac{A_1'B}{A_1'C}.$$

In view of *1$^{st}$ Proposition*, we get that:

$$\frac{A_1'B}{A_1'C} = \left(\frac{AB}{AC}\right)^2.$$

# 7<sup>th</sup> Proposition.

The tangents taken to the extremes of a diagonal of a circle circumscribed to the harmonic quadrilateral intersect on the other diagonal.

## *Proof.*

Let $P$ be the intersection of a tangent taken in $D$ to the circle circumscribed to the harmonic quadrilateral $ABCD$ with $AC$ (see *Figure* 4).

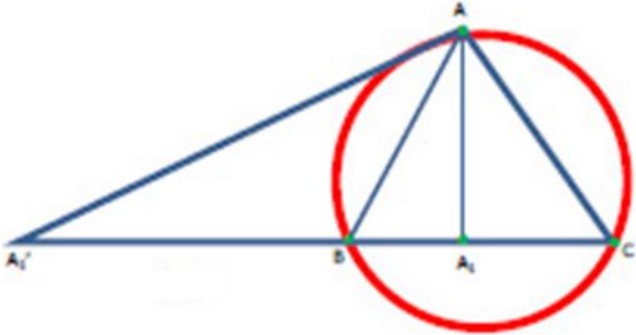

*Figure 4.*

Since triangles PDC and PAD are alike, we conclude that:

$$\frac{PD}{PA} = \frac{PC}{PD} = \frac{DC}{AD} \ . \tag{5}$$

From relations (5), we find that:

$$\frac{PA}{PC} = \left(\frac{AD}{DC}\right)^2 . \tag{6}$$

This relationship indicates that P is the harmonic conjugate of $K$ with respect to $A$ and $C$, so $DP$ is an external simedian from $D$ of the triangle $ADC$.

Similarly, if we denote by $P'$ the intersection of the tangent taken in $B$ to the circle circumscribed with $AC$, we get:

$$\frac{P'A}{P'C} = \left(\frac{BA}{BC}\right)^2. \tag{7}$$

From (6) and (7), as well as from the properties of the harmonic quadrilateral, we know that:

$$\frac{AB}{BC} = \frac{AD}{DC},$$

which means that:

$$\frac{PA}{PC} = \frac{P'A}{P'C},$$

hence $P = P'$.

Similarly, it is shown that the tangents taken to $A$ and $C$ intersect at point $Q$ located on the diagonal $BD$.

## $5^{th}$ *Remark.*

a. The points $P$ and $Q$ are the diagonal poles of $BD$ and $AC$ in relation to the circle circumscribed to the quadrilateral.

b. From the previous *Proposition*, it follows that in a triangle the internal simedian of an angle is consecutive to the external simedians of the other two angles.

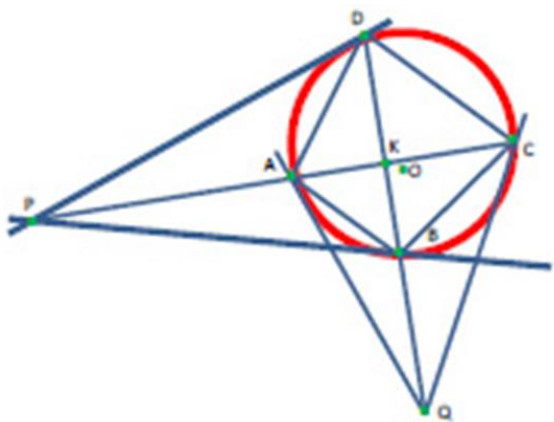

*Figure 5.*

# 8<sup>th</sup> Proposition.

Let $ABCD$ be an harmonic quadrilateral inscribed in the circle of center $O$ and let $P$ and $Q$ be the intersections of the tangents taken in $B$ and $D$, respectively in $A$ and $C$ to the circle circumscribed to the quadrilateral.

If $\{K\} = AC \cap BD$, then the orthocenter of triangle $PKQ$ is $O$.

### *Proof.*

From the properties of tangents taken from a point to a circle, we conclude that $PO \perp BD$ and $QO \perp AC$. These relations show that in the triangle $PKQ$, $PO$ and $QO$ are heights, so $O$ is the orthocenter of this triangle.

# 4ᵗʰ Definition.

The Apollonius's circle related to the vertex $A$ of the triangle $ABC$ is the circle built on the segment $[DE]$ in diameter, where $D$ and $E$ are the feet of the internal, respectively external, bisectors taken from $A$ to the triangle $ABC$.

## 6ᵗʰ Remark.

If the triangle $ABC$ is isosceles with $AB = AC$, the Apollonius's circle corresponding to vertex $A$ is not defined.

# 9ᵗʰ Proposition.

The Apollonius's circle relative to the vertex $A$ of the triangle $ABC$ has as center the feet of the external simedian taken from $A$.

## Proof.

Let $O_a$ be the intersection of the external simedian of the triangle $ABC$ with $BC$ (see *Figure 6*).

Assuming that $m(\hat{B}) > m(\hat{C})$, we find that:

$m(\widehat{EAB}) = \frac{1}{2}[m(\hat{B}) + m(\hat{C})]$.

$O_a$ being a tangent, we find that $m(\widehat{O_aAB}) = m(\hat{C})$.

Withal,

$m(EAO_a) = \frac{1}{2}[m(\hat{B}) - m(\hat{C})]$

and

$$m(AEO_a) = \frac{1}{2}\big[m(\hat{B}) - m(\hat{C})\big].$$

It results that $O_aE = O_aA$; onward, $EAD$ being a right angled triangle, we obtain: $O_aA = O_aD$; hence $O_a$ is the center of Apollonius's circle corresponding to the vertex $A$.

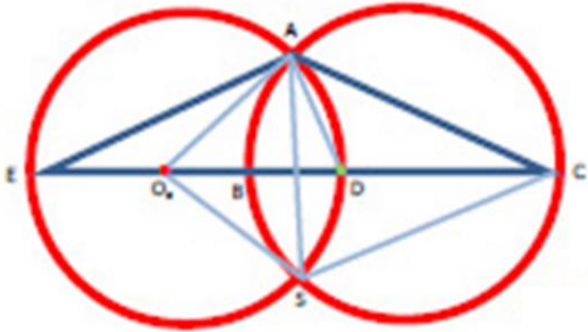

*Figura 6.*

## 10$^{th}$ Proposition.

Apollonius's circle relative to the vertex $A$ of triangle $ABC$ cuts the circle circumscribed to the triangle following the internal simedian taken from $A$.

### *Proof.*

Let $S$ be the second point of intersection of Apollonius's Circle relative to vertex $A$ and the circle circumscribing the triangle $ABC$.

Because $O_a A$ is tangent to the circle circumscribed in A, it results, for reasons of symmetry, that $O_a S$ will be tangent in $S$ to the circumscribed circle.

For triangle $ACS$, $O_a A$ and $O_a S$ are external simedians; it results that $CO_a$ is internal simedian in the triangle $ACS$, furthermore, it results that the quadrilateral $ABSC$ is an harmonic quadrilateral.

Consequently, $AS$ is the internal simedian of the triangle $ABC$ and the property is proven.

### $7^{th}$ Remark.

From this, in view of *Figure 5*, it results that the circle of center $Q$ passing through $A$ and $C$ is an Apollonius's circle relative to the vertex $A$ for the triangle $ABD$. This circle (of center $Q$ and radius $QC$) is also an Apollonius's circle relative to the vertex $C$ of the triangle $BCD$.

Similarly, the Apollonius's circle corresponding to vertexes $B$ and $D$ and to the triangles ABC, and ADC respectively, coincide.

We can now formulate the following:

## 11$^{th}$ Proposition.

In an harmonic quadrilateral, the Apollonius's circle - associated with the vertexes of a diagonal and to the triangles determined by those vertexes to the other diagonal - coincide.

Radical axis of the Apollonius's circle is the right determined by the center of the circle circumscribed to the harmonic quadrilateral and by the intersection of its diagonals.

*Proof.*

Referring to *Fig. 5*, we observe that the power of *O* towards the Apollonius's Circle relative to vertexes *B* and *C* of triangles *ABC* and *BCU* is:

$OB^2 = OC^2$.

So *O* belongs to the radical axis of the circles.

We also have $KA \cdot KC = KB \cdot KD$, relatives indicating that the point *K* has equal powers towards the highlighted Apollonius's circle.

# References.

# Triangulation of a Triangle with Triangles having Equal Inscribed Circles

In this article, we solve the following problem:

*Any triangle can be divided by a cevian in two triangles that have congruent inscribed circles.*

## Solution.

We consider a given triangle $ABC$ and we show that there is a point $D$ on the side $(BC)$ so that the inscribed circles in the triangles $ABD$ , $ACD$ are congruent. If $ABC$ is an isosceles triangle $(AB = AC)$, where $D$ is the middle of the base $(BC)$, we assume that $ABC$ is a non-isosceles triangle. We note $I_1, I_2$ the centers of the inscribed congruent circles; obviously, $I_1 I_2$ is parallel to the $BC$. (1)

We observe that:

$$m(\sphericalangle I_1 A I_2) = \frac{1}{2} m(\hat{A}). \qquad (2)$$

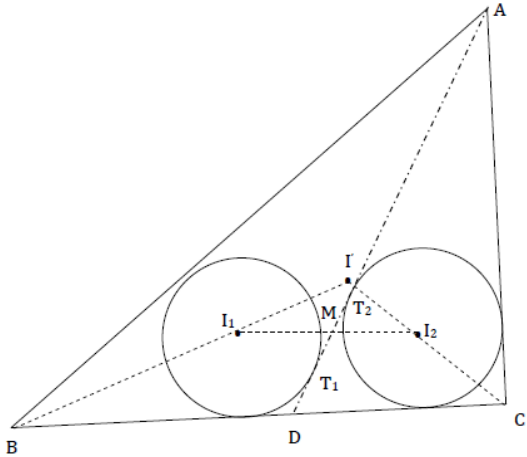

*Figure 1.*

If $T_1, T_2$ are contacts with the circles of the cevian $AD$ , we have $\Delta\, I_1 T_1 M \equiv \Delta\, I_2 T_2 M$ ; let $M$ be the intersection of $I_1 I_2$ with $AD$, see *Figure 1*.

From this congruence, it is obvious that:

$(I_1 M) \equiv (I_2 M).$ (3)

Let $I$ be the center of the circle inscribed in the triangle $ABC$; we prove that: $AI$ is a simedian in the triangle $I_1 A I_2$. (4)

Indeed, noting $\alpha = m\big(\widehat{BAI_1}\big)$ , it follows that $m(\sphericalangle I_1 AM) = \alpha$. From $\sphericalangle I_1 A I_2 = \sphericalangle BAI,$ it follows that $\sphericalangle BAI_1 \equiv \sphericalangle IAI_2$ , therefore $\sphericalangle I_1 AM \equiv \sphericalangle IAI_2$ , indicating that $AM$ and $AI$ are isogonal cevians in the triangle $I_1 A I_2$.

Since in this triangle $AM$ is a median, it follows that $AI$ is a bimedian.

Now, we show how we build the point $D$, using the conditions (1) – (4), and then we prove that this construction satisfies the enunciation requirements.

## Building the point $D$.

1°: We build the circumscribed circle of the given triangle $ABC$; we build the bisector of the angle $BAC$ and denote by P its intersection with the circumscribed circle (see *Figure 2*).

2°: We build the perpendicular on $C$ to $CP$ and $(BC)$ side mediator; we denote $O_1$ the intersection of these lines.

3°: We build the circle $C(O_1; O_1C)$ and denote $A'$ the intersection of this circle with the bisector $AI$ ($A'$ is on the same side of the line $BC$ as $A$).

4°: We build through $A$ the parallel to $A'O_1$ and we denote it $IO_1$.

5°: We build the circle $C(O_1'; O_1'A)$ and we denote $I_1, I_2$ its intersections with $BI$, and $CI$ respectively.

6°: We build the middle $M$ of the segment $(I_1I_2)$ and denote by $D$ the intersection of the lines $AM$ and $BC$.

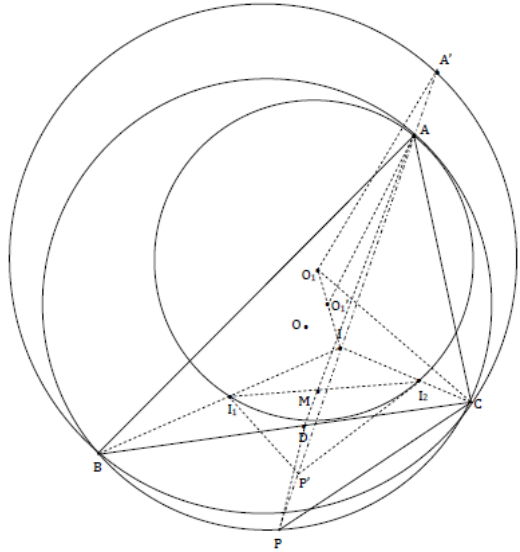

*Figure 2.*

## Proof.

The point $P$ is the middle of the arc $\overline{BC}$, then $m(\widehat{PCB}) = \frac{1}{2}m(\hat{A})$.

The circle $\mathsf{C}(O_1; O_1C)$ contains the arc from which points the segment (BC) „is shown" under angle measurement $\frac{1}{2}m(\hat{A})$.

The circle $\mathsf{C}(O_1'; O_1'A)$ is homothetical to the circle $\mathsf{C}(O_1; O_1C)$ by the homothety of center $I$ and by the report $\frac{IA'}{IA}$; therefore, it follows that $I_1I_2$ will be parallel to the $BC$, and from the points of circle $\mathsf{C}(O_1'; O_1'A)$ of the same side of $BC$ as $A$, the segment $(I_1 I_2)$ „is shown" at an angle of measure $\frac{1}{2} m(\hat{A})$. Since the tangents taken in $B$ and $C$ to the circle $C(O_1; O_1 C)$ intersect in $P$, on the bisector $AI$, as from a property of simedians, we get that $A'I$ is a simedian in the triangle $A'BC$.

Due to the homothetical properties, it follows also that the tangents in the points $I_1, I_2$ to the circle $C(O_1'; O_1'A)$ intersect in a point $P'$ located on $AI$, *i.e.* $AP'$ contains the simedian $(AS)$ of the triangle $I_1 A I_2$, noted $\{S\} = AP' \cap I_1 I_2$.

In the triangle $I_1 A I_2$, $AM$ is a median, and $AS$ is simedian, therefore $\sphericalangle I_1 AM \equiv I_2 AS$; on the other hand, $\sphericalangle BAS \equiv \sphericalangle I_1 A I_2$; it follows that $\sphericalangle BAI_1 \equiv I_2 AS$, and more: $\sphericalangle I_1 AM \equiv BAI_1$, which shows that $AI_1$ is a bisector in the triangle $BAD$; plus, $I_1$, being located on the bisector of the angle $B$, it follows that this point is the center of the circle inscribed in the triangle $BAD$.

Analogous considerations lead to the conclusion that $I_2$ is the center of the circle inscribed in the triangle $ACD$. Because $I_1 I_2$ is parallel to $BC$, it follows that the rays of the circles inscribed in the triangles $ABD$ and $ACD$ are equal.

### Discussion.

The circles $C(O_1; O_1 A')$, $C(O_1'; O_1'A)$ are unique; also, the triangle $I_1 A I_2$ is unique; therefore, determined as before, the point $D$ is unique.

*Remark.*

At the beginning of the *Proof*, we assumed that $ABC$ is a non-isosceles triangle with the stated property. There exist such triangles; we can construct such a triangle starting "backwards". We consider two given congruent external circles and, by tangent constructions, we highlight the $ABC$ triangle.

## Open problem.

Given a scalene triangle $ABC$ , could it be triangulated by the cevians $AD$ , $AE$ , with $D$ , $E$ belonging to $(BC)$, so that the inscribed circles in the triangles $ABD$, $DAE$ and the $EAC$ to be congruent?

# An Application of a Theorem of Orthology

In this short article, we make connections between Problem 21 of [1] and the theory of orthological triangles.

The enunciation of the competitional problem is:

Let $(T_A)$, $(T_B)$, $(T_C)$ be the tangents in the peaks $A, B, C$ of the triangle $ABC$ to the circle circumscribed to the triangle. Prove that the perpendiculars drawn from the middles of the opposite sides on $(T_A)$, $(T_B)$, $(T_C)$ are concurrent and determine their concurrent point.

We formulate and we demonstrate below a sentence containing in its proof the solution of the competitional problem in this case.

## Proposition.

The tangential triangle and the median triangle of a given triangle $ABC$ are orthological. The orthological centers are $O$ – the center of the circle circumscribed to the triangle $ABC$, and $O_9$ – the center of the circle of $ABC$ triangle's 9 points.

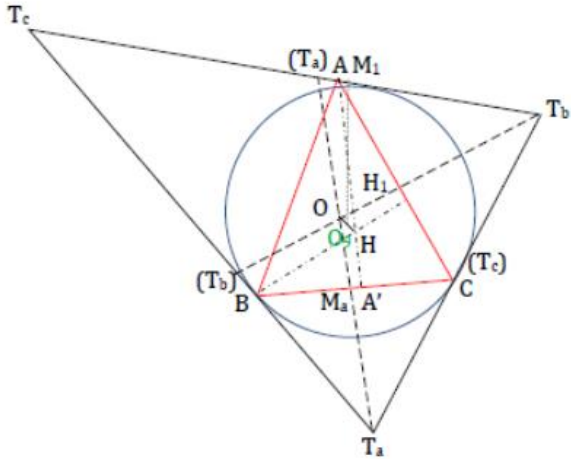

*Figure 1.*

## Proof.

Let $M_a M_b M_c$ be the median triangle of triangle $ABC$ and $T_a T_b T_c$ the tangential triangle of the triangle $ABC$. It is obvious that the triangle $T_a T_b T_c$ and the triangle $ABC$ are orthological and that $O$ is the orthological center.

Verily, the perpendiculars taken from $T_a, T_b, T_c$ on $BC$; $CA$; $AB$ respectively are internal bisectors in the triangle $T_a T_b T_c$ and consequently passing through $O$, which is center of the circle inscribed in triangle $T_a T_b T_c$. Moreover, $T_a O$ is the mediator of $(BC)$ and accordingly passing through $M_a$, and $T_a M_c$ is perpendicular on $BC$, being a mediator, but also on $M_b M_c$ which is a median line.

From orthological triangles theorem, it follows that the perpendiculars taken from $M_a, M_b, M_c$ on $T_bT_c$, $T_cT_a$, $T_aT_b$ respectively, are concurrent. The point of concurrency is the second orthological center of triangles $T_aT_bT_c$ and $M_aM_bM_c$. We prove that this point is $O_9$ – the center of Euler circle of triangle $ABC$. We take $M_aM_1 \perp T_bT_c$ and denote by $\{H_1\} = M_aM_1 \cap AH$, $H$ being the orthocenter of the triangle $ABC$. We know that $AH = 2OM_a$; we prove this relation vectorially.

From Sylvester's relation, we have that: $\overrightarrow{OH} = \overrightarrow{OA} + \overrightarrow{OB} + \overrightarrow{OC}$, but $\overrightarrow{OB} + \overrightarrow{OC} = 2\overrightarrow{OM_a}$; it follows that $\overrightarrow{OH} - \overrightarrow{OA} = 2\overrightarrow{OM_a}$, so $\overrightarrow{AH} = 2\overrightarrow{OM_a}$; changing to module, we have $AH = 2OM_a$. Uniting $O$ to $A$, we have $OA \perp T_bT_c$, and because $M_aM_1 \perp T_bT_c$ and $AH \parallel OM_a$, it follows that the quadrilateral $OM_aH_1A$ is a parallelogram.

From $AH_1 = OM_a$ and $AH = 2OM_a$ we get that $H_1$ is the middle of $(AH)$, so $H_1$ is situated on the circle of the 9 points of triangle $ABC$. On this circle, we find as well the points $A'$ - the height foot from $A$ and $M_a$; since $\sphericalangle AA'M_a = 90^0$, it follows that $M_aH_1$ is the diameter of Euler circle, therefore the middle of $(M_aH_1)$, which we denote by $O_9$, is the center of Euler's circle; we observe that the quadrilateral $H_1HM_aO$ is as well a parallelogram; it follows that $O_9$ is the middle of segment $[OH]$.

In conclusion, the perpendicular from $M_a$ on $T_bT_c$ pass through $O_9$.

Analogously, we show that the perpendicular taken from $M_b$ on $T_a T_c$ pass through $O_9$ and consequently $O_9$ is the enunciated orthological center.

### *Remark.*

The triangles $M_a M_b M_c$ and $T_a T_b T_c$ are homological as well, since $T_a M_a$ , $T_b M_b$ , $T_c M_c$ are concurrent in O, as we already observed, therefore the triangles $T_a T_b T_c$ and $M_a M_b M_c$ are orthohomological of rank I (see [2]).

From P. Sondat theorem (see [4]), it follows that the Euler line $O O_9$ is perpendicular on the homological axis of the median triangle and of the tangential triangle corresponding to the given triangle $ABC$.

### *Note.*

(Regarding the triangles that are simultaneously orthological and homological)

In the article *A Theorem about Simultaneous Orthological and Homological Triangles*, by Ion Patrascu and Florentin Smarandache, we stated and proved a theorem which was called by Mihai Dicu in [2] and [3] *The Smarandache-Patrascu Theorem of Orthohomological Triangle*; then, in [4], we used the term *ortohomological triangles* for the triangles that are simultaneously orthological and homological.

The term ortohomological triangles was used by J. Neuberg in *Nouvelles Annales de Mathematiques*

(1885) to name the triangles that are simultaneously orthogonal (one has the sides perpendicular to the sides of the other) and homological.

We suggest that the triangles that are simultaneously orthogonal and homological to be called ortohomological triangles of first rank, and triangles that are simultaneously orthological and homological to be called ortohomological triangles of second rank.

# References.

# The Dual of a Theorem Relative to the Orthocenter of a Triangle

In [1] we introduced the notion of **Bobillier's transversal relative to a point *O* in the plane of a triangle *ABC*** ; we use this notion in what follows.

We transform by duality with respect to a circle $\mathbb{C}(o,r)$ the following theorem relative to the orthocenter of a triangle.

## 1$^{\text{st}}$ Theorem.

If $ABC$ is a nonisosceles triangle, $H$ its orthocenter, and $AA_1, BB_1, CC_1$ are cevians of a triangle concurrent at point $Q$ different from $H$, and $M, N, P$ are the intersections of the perpendiculars taken from $H$ on given cevians respectively, with $BC, CA, AB$, then the points $M, N, P$ are collinear.

*Proof.*

We note with $\alpha = m(\sphericalangle BAA_1)$; $\beta = m(\sphericalangle CBB_1), \gamma = m(\sphericalangle ACC_1)$, see *Figure 1*.

According to Ceva's theorem, trigonometric form, we have the relation:

$$\frac{\sin\alpha}{\sin(A-\alpha)} \cdot \frac{\sin\beta}{\sin(B-\beta)} \cdot \frac{\sin\gamma}{\sin(C-\gamma)} = 1. \tag{1}$$

We notice that:

$$\frac{MB}{MC} = \frac{\text{Area}(MHB)}{\text{Area}(MHC)} = \frac{MH \cdot HB \cdot \sin\sin(\widehat{MHB})}{MH \cdot HC \cdot \sin(\widehat{MHC})} \, .$$

Because:

$$\sphericalangle MHB \equiv \sphericalangle A_1AC,$$

as angles of perpendicular sides, it follows that

$$m(\sphericalangle MHB) = m(\hat{A}) - \alpha.$$

Therewith:

$$m(\sphericalangle MHC) = m(\widehat{MHB}) + m(\widehat{BHC}) = 180^0\alpha.$$

We thus get that:

$$\frac{MB}{MC} = \frac{\sin(A-\alpha)}{\sin\alpha} \cdot \frac{HB}{HC}.$$

Analogously, we find that:

$$\frac{NC}{NA} = \frac{\sin(B-\beta)}{\sin\beta} \cdot \frac{HC}{HA} \, ;$$

$$\frac{PA}{PB} = \frac{\sin(C-\gamma)}{\sin\gamma} \cdot \frac{HA}{HB} \, .$$

Applying the reciprocal of Menelaus' theorem, we find, in view of (1), that:

$$\frac{MB}{MC} \cdot \frac{HC}{HA} \cdot \frac{PA}{PB} = 1.$$

This shows that $M, N, P$ are collinear.

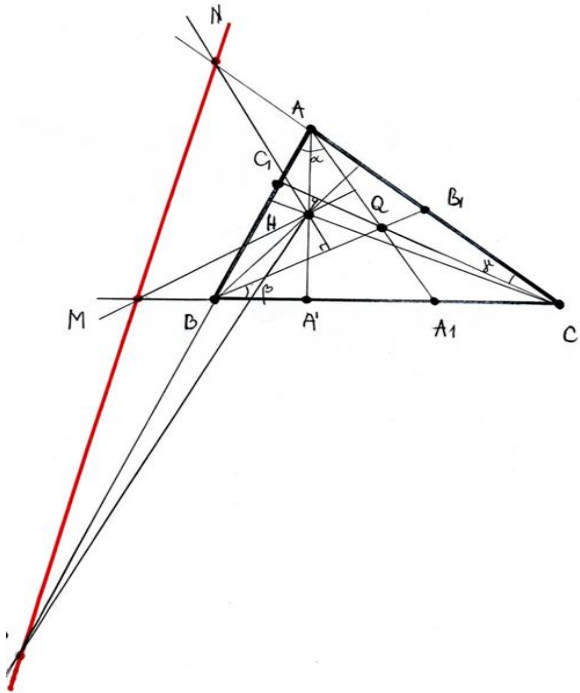

*Figure 1.*

*Note.*

$1^{st}$ *Theorem* is true even if $ABC$ is an obtuse, nonisosceles triangle.

The proof is adapted analogously.

## $2^{nd}$ **Theorem.**

(The Dual of the Theorem 1)

If $ABC$ is a triangle, $O$ a certain point in his plan, and $A_1, B_1, C_1$ Bobillier's transversals relative to $O$ of $ABC$ triangle, as well as $A_2 - B_2 - C_2$ a certain transversal in $ABC$, and the perpendiculars in $O$, and on $OA_2, OB_2, OC_2$ respectively, intersect the Bobillier's transversals in the points $A_3, B_3, C_3$, then the ceviens $AA_3, BB_3, CC_3$ are concurrent.

## *Proof.*

We convert by duality with respect to a circle $\mathbb{C}(o, r)$ the figure indicated by the statement of this theorem, i.e. *Figure 2*. Let $a, b, c$ be the polars of the points $A, B, C$ with respect to the circle $\mathbb{C}(o, r)$. To the lines $BC, CA, AB$ will correspond their poles $\{A'4 = bnc; \{B'4 = cna; \{C'4 = anb.$

To the points $A_1, B_1, C_1$ will respectively correspond their polars $a_1, b_1, c_1$ concurrent in transversal's pole $A_1 - B_1 - C_1$.

Since $OA_1 \perp OA$, it means that the polars $a$ and $a_1$ are perpendicular, so $a_1 \perp B'C'$, but $a_1$ pass through $A'$, which means that $Q'$ contains the height from $A'$ of $A'B'C'$ triangle and similarly $b_1$ contains the height from $B'$ and $c_1$ contains the height from $C'$ of $A'B'C'$ triangle.

Consequently, the pole of $A_1 - B_1 - C_1$ transversal is the orthocenter $H'$ of $A'B'C'$ triangle. In the same way, to the points $A_2, B_2, C_2$ will correspond the polars to $a_2, b_2, c_2$ which pass respectively through $A', B', C'$ and are concurrent in a point $Q'$, the pole of the line $A_2 - B_2 - C_2$ with respect to the circle $\mathbb{C}(o, r)$.

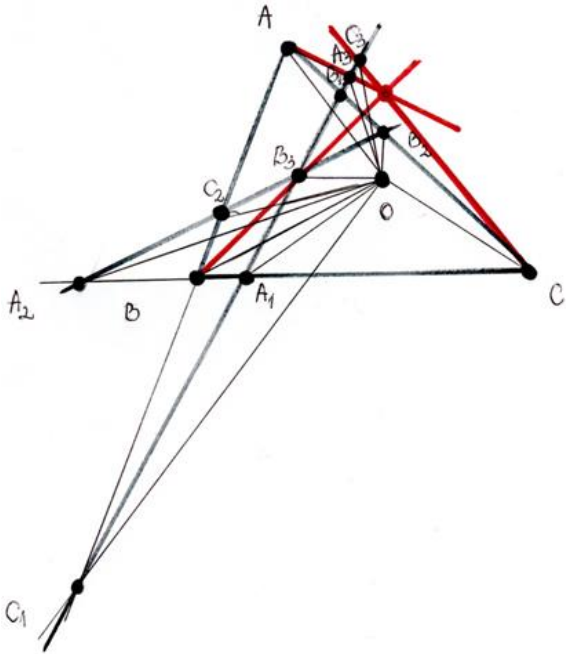

Given $OA_2 \perp OA_3$, it means that the polar $a_2$ and $a_3$ are perpendicular, $a_2$ correspond to the cevian $A'Q'$, also $a_3$ passes through the the pole of the transversal $A_1 - B_1 - C_1$, so through $H'$, in other words $Q_3$ is perpendicular taken from $H'$ on $A'Q'$; similarly, $b_2 \perp b_3, c_2 \perp c_3$, so $b_3$ is perpendicular taken from $H'$ on $C'Q'$. To the cevian $AA_3$ will correspond by duality considered to its pole, which is the intersection of the polars of $A$ and $A_3$, i.e. the intersection of lines $a$ and $a_3$ , namely the intersection of $B'C'$ with the perpendicular taken from $H'$ on $A'Q'$; we denote this point by $M'$. Analogously, we get the points $N'$ and $P'$. Thereby, we got the configuration from *1st Theorem* and *Figure 1*, written for triangle $A'B'C'$ of orthocenter $H'$. Since from *1st Theorem* we know that $M', N', P'$ are collinear, we get the the cevians $AA_3, BB_3, CC_3$ are concurrent in the pole of transversal $M' - N' - P'$ with respect to the circle $\mathbb{C}(o, r)$, and *2nd Theorem* is proved.

# References.

# The Dual Theorem Concerning Aubert's Line

In this article we introduce the concept of **Bobillier's transversal of a triangle with respect to a point in its plan**; we prove the Aubert's Theorem about the collinearity of the orthocenters in the triangles determined by the sides and the diagonals of a complete quadrilateral, and we obtain the Dual Theorem of this Theorem.

## 1$^{st}$ Theorem.

(E. Bobillier)

Let $ABC$ be a triangle and $M$ a point in the plane of the triangle so that the perpendiculars taken in $M$, and $MA, MB, MC$ respectively, intersect the sides $BC, CA$ and $AB$ at $Am$, $Bm$ şi $Cm$. Then the points $Am$, $Bm$ and $Cm$ are collinear.

*Proof.*

We note that $\frac{AmB}{AmC} = \frac{\text{aria}\,(BMAm)}{\text{aria}\,(CMAm)}$ (see *Figure 1*).

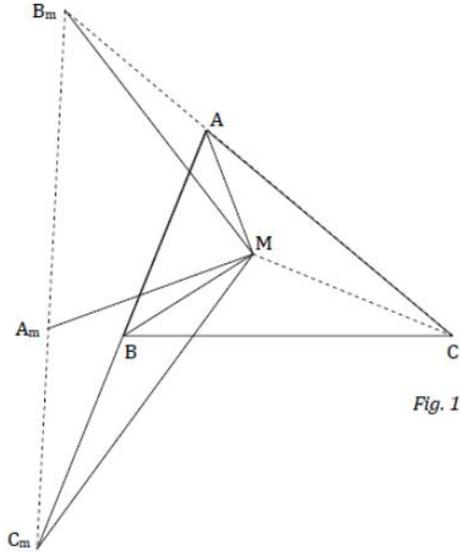

*Figure 1.*

Area $(BMAm) = \frac{1}{2} \cdot BM \cdot MAm \cdot \sin(B\widehat{MA}m)$.

Area $(CMAm) = \frac{1}{2} \cdot CM \cdot MAm \cdot \sin(C\widehat{MA}m)$.

Since

$$m(C\widehat{MA}m) = \frac{3\pi}{2} - m(\widehat{AMC}),$$

it explains that

$$\sin(C\widehat{MA}m) = -\cos(\widehat{AMC});$$

$$\sin(B\widehat{MA}m) = \sin\left(\widehat{AMB} - \frac{\pi}{2}\right) = -\cos(\widehat{AMB}).$$

Therefore:

$$\frac{AmB}{AmC} = \frac{MB \cdot \cos(\widehat{AMB})}{MC \cdot \cos(\widehat{AMC})} \ (1).$$

In the same way, we find that:

$$\frac{BmC}{BmA} = \frac{MC}{MA} \cdot \frac{\cos(\widehat{BMC})}{\cos(\widehat{AMB})} \ (2);$$

$$\frac{CmA}{CmB} = \frac{MA}{MB} \cdot \frac{\cos(\widehat{AMC})}{\cos(\widehat{BMC})} \ (3).$$

The relations (1), (2), (3), and the reciprocal Theorem of Menelaus lead to the collinearity of points $Am, Bm, Cm$.

*Note.*

Bobillier's Theorem can be obtained – by converting the duality with respect to a circle – from the theorem relative to the concurrency of the heights of a triangle.

# 1ˢᵗ Definition.

It is called Bobillier's transversal of a triangle $ABC$ with respect to the point $M$ the line containing the intersections of the perpendiculars taken in $M$ on $AM$, $BM$, and $CM$ respectively, with sides $BC$, CA and $AB$.

*Note.*

The Bobillier's transversal is not defined for any point $M$ in the plane of the triangle $ABC$, for example, where $M$ is one of the vertices or the orthocenter $H$ of the triangle.

## 2nd Definition.

If $ABCD$ is a convex quadrilateral and $E, F$ are the intersections of the lines $AB$ and $CD$, $BC$ and $AD$ respectively, we say that the figure $ABCDEF$ is a complete quadrilateral. The complete quadrilateral sides are $AB, BC, CD, DA$,  and $AC, BD$ and $EF$ are diagonals.

## 2nd Theorem.

(Newton-Gauss)

The diagonals' middles of a complete quadrilateral are three collinear points. To prove 2nd theorem, refer to [1].

### *Note.*

It is called *Newton-Gauss Line* of a quadrilateral the line to which the diagonals' middles of a complete quadrilateral belong.

## 3rd Theorem.

(Aubert)

If $ABCDEF$ is a complete quadrilateral, then the orthocenters $H_1, H_2, H_3, H_4$ of the traingles $ABF$, $AED$, $BCE$, and $CDF$ respectively, are collinear points.

*Proof.*

Let $A_1, B_1, F_1$ be the feet of the heights of the triangle $ABF$ and $H_1$ its orthocenter (see *Fig. 2*). Considering the power of the point $H_1$ relative to the circle circumscribed to the triangle ABF, and given the triangle orthocenter's property according to which its symmetrics to the triangle sides belong to the circumscribed circle, we find that:

$H_1A \cdot H_1A_1 = H_1B \cdot H_1B_1 = H_1F \cdot H_1F_1.$

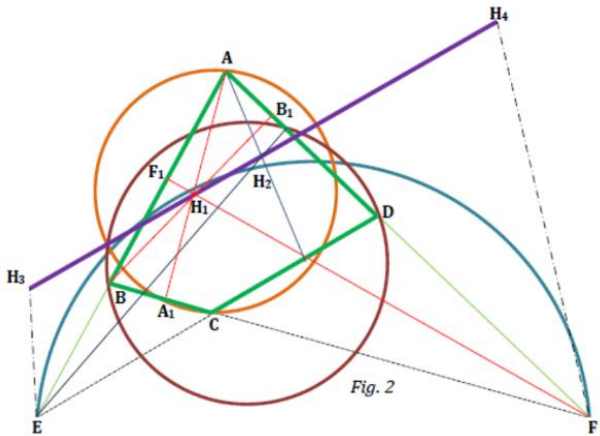

*Figure 2.*

This relationship shows that the orthocenter $H_1$ has equal power with respect to the circles of diameters $[AC], [BD], [EF]$. As well, we establish that the orthocenters $H_2, H_3, H_4$ have equal powers to these circles. Since the circles of diameters $[AC], [BD], [EF]$ have collinear centers (belonging to the Newton-Gauss line of the $ABCDEF$ quadrilateral), it follows that the points $H_1, H_2, H_3, H_4$ belong to the radical axis of the circles, and they are, therefore, collinear points.

*Notes.*

1. It is called *the Aubert Line* or the line of the complete quadrilateral's orthocenters the line to which the orthocenters $H_1, H_2, H_3, H_4$ belong.

2. The Aubert Line is perpendicular on the Newton-Gauss line of the quadrilateral (radical axis of two circles is perpendicular to their centers' line).

## $4^{\text{th}}$ Theorem.

(The Dual Theorem of the $3^{\text{rd}}$ Theorem)

If $ABCD$ is a convex quadrilateral and $M$ is a point in its plane for which there are the Bobillier's transversals of triangles $ABC$, $BCD$, $CDA$ and $DAB$; thereupon these transversals are concurrent.

*Proof.*

Let us transform the configuration in *Fig. 2,* by duality with respect to a circle of center $M$.

By the considered duality, the lines $a, b, c, d, e$ şi $f$ correspond to the points $A, B, C, D, E, F$ (their polars).

It is known that polars of collinear points are concurrent lines, therefore we have: $a \cap b \cap e = \{A'\}$,

$b \cap c \cap f = \{B'\}$, $c \cap d \cap e = \{C'\}$, $d \cap f \cap a = \{D'\}$, $a \cap c = \{E'\}$, $b \cap d = \{F'\}$.

Consequently, by applicable duality, the points $A', B', C', D', E'$ and $F'$ correspond to the straight lines $AB, BC, CD, DA, AC, BD$.

To the orthocenter $H_1$ of the triangle $AED$, it corresponds, by duality, its polar, which we denote $A'_1 - B'_1 - C'_1$, and which is the Bobillier's transversal of the triangle $A'C'D'$ in relation to the point $M$. Indeed, the point $C'$ corresponds to the line $ED$ by duality; to the height from $A$ of the triangle $AED$, also by duality, it corresponds its pole, which is the point $C'_1$ located on $A'D'$ such that $m(\widehat{C'MC'_1}) = 90^0$.

To the height from $E$ of the triangle $AED$, it corresponds the point $B'_1 \in A'C'$ such that $m(\widehat{D'MB'_1}) = 90^0$.

Also, to the height from $D$, it corresponds $A'_1 \in C'D'$ on C such that $m(A'MA'_1) = 90^0$. To the orthocenter $H_2$ of the triangle $ABF$, it will correspond, by applicable duality, the Bobillier's transversal $A'_2 - B'_2 - C'_2$ in the triangle $A'B'D'$ relative to the point $M$. To the orthocenter $H_3$ of the triangle $BCE$, it will correspond the Bobillier's transversal $A_3' - B'_3 - C_3'$ in the triangle $A'B'C'$ relative to the point $M$, and to the orthocenter $H_4$ of the triangle $CDF$, it will correspond the transversal $A'_4 - B'_4 - C'_4$ in the triangle $C'D'B'$ relative to the point $M$.

The Bobillier's transversals $A_i' - B_i' - C_i'$, $i = \overline{1,4}$ correspond to the collinear points $H_i$, $i = \overline{1,4}$.

These transversals are concurrent in the pole of the line of the orthocenters towards the considered duality.

It results that, given the quadrilateral $A'B'C'D'$, the Bobillier's transversals of the triangles $A'C'D'$, $A'B'D'$, $A'B'C'$ and $C'D'B'$ relative to the point $M$ are concurrent.

# References.

We approach several themes of classical geometry of the circle and complete them with some original results, showing that not everything in traditional math is revealed, and that it still has an open character.

The topics were chosen according to authors' aspiration and attraction, as a poet writes lyrics about spring according to his emotions.

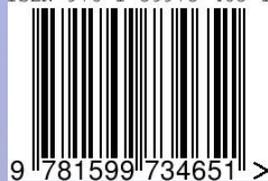